\documentclass{amsart}
\usepackage[english]{babel}
\usepackage[latin1]{inputenc}
\usepackage[dvips,final]{graphics}
\usepackage{amsmath,amsfonts,amssymb,amsthm,amscd,array}
\usepackage[vcentermath]{youngtab}
\usepackage{young}
\usepackage{pstricks}
\usepackage{textcomp}
 \usepackage[final]{epsfig}
\theoremstyle{plain}
\newtheorem{thm}{Theorem}[section]
\newtheorem{cor}[thm]{Corollary}
\newtheorem{lem}[thm]{Lemma}
\newtheorem{prop}[thm]{Proposition}

\theoremstyle{definition}
\newtheorem{defn}[thm]{Definition}
\newtheorem{rem}[thm]{Remark}
\newtheorem{ex}[thm]{Example}

\newcommand{\R}{\mathbb{R}}
\newcommand{\Z}{\mathbb{Z}}
\newcommand{\C}{\mathbb{C}}
\newcommand{\N}{\mathbb{N}}

\newcommand{\Q}{\mathbb{Q}}

\newcommand{\B}{\mathcal{B}}
\newcommand{\F}{\mathcal{F}}
\newcommand{\cc}{\mathcal{C}}
\newcommand{\f}{\mathcal{L}} 

\newcommand{\gn}{\mathfrak{n}}
\newcommand{\gog}{\mathfrak{g}}

\newcommand{\gb}{\mathfrak{b}}
\newcommand{\gh}{\mathfrak{h}}

\newcommand{\U}{\mathcal{U}}

\newcommand{\ii}{\textup{\bf{i}}}
\newcommand{\tr}{\textup{Trop }}
\newcommand{\Gr}{\textup{Gr}}

\begin{document}

\title[Geometric lifting of the canonical basis ...]{Geometric lifting of the canonical basis and semitoric degenerations of Richardson varieties}

\author{Sophie Morier-Genoud}

\address{D\'epartement de Math\'ematiques, Universit\'e Claude Bernard Lyon I,
69622 Villeurbanne Cedex, France}

\email{morier@math.univ-lyon1.fr}

\date{April 26, 2005.}

\subjclass{14M25, 16W35, 14M15}

\begin{abstract}
In the $\mathfrak{sl}_n$ case, A. Berenstein and A. Zelevinsky studied in \cite{berensteinzelevinsky1} the
Sch\"utzenberger involution in terms of Lusztig's canonical basis. We
generalize their construction and formulas for any semisimple Lie algebra.
We use the geometric lifting of the canonical basis, on which an
analogue of the Sch\"utzenberger involution can be given. As an application,
we construct semitoric degenerations of Richardson varieties, following a method of P.~Caldero,~\cite{Philippe}.
\end{abstract}

\maketitle

\vspace{-0.8cm}
\tableofcontents

\section*{Introduction}
Let $G$ be a semisimple simply connected complex Lie group. Fix opposite Borel subgroups $B$ and $B^-$ of $G$. In this paper we consider subvarieties of the flag variety $G/B$ known as Richardson varieties.
They first appear in \cite{richardson}. Our problem is to construct toric or semitoric degenerations of these varieties. Such constructions have already be done in the special cases of the flag variety and the Schubert varieties, see \cite{GL}, \cite{chirivi} and \cite{Philippe}. Our approach consists to extend the method introduced by \cite{Philippe}. Let us mention that the method of \cite{Philippe} was recently extend for the degenerations of spherical varieties, see \cite{alexeevbrion}.\par
A Richardson variety $X_w^\tau$ is the intersection of a Schubert variety $X_w:=\overline {BwB/B}$ and an opposite Schubert variety $X^{-}_{\tau }:=\overline {B^-{\tau }B/B}$, where $w$ and $\tau$ are elements in the Weyl group $W$ of $G$. The opposite Schubert variety $X^\tau$ is the image of a Schubert variety under the action of the longest element $w_0$ of $W$. This element plays an important role in our study. In order to construct toric degenerations of these varieties, we define filtrations on the homogenous coordinates algebras associated to the varieties (Sections \ref{SectionG/B}, \ref{SectionXw} and \ref{SectionXwtau}). All these algebras are direct sums of subspaces of $G$-modules. The algebra $R^\tau$ associated to the opposite Schubert variety $X^\tau$ is related to the algebra $R_\tau$ associated to the Schubert variety $X_\tau$ via the action of $w_0$. It is then important to understand the action of $w_0$ on the $G$-modules.\par 
An important tool in our work is the canonical/global basis of Lusztig and Kashiwara. This basis $\B$ lays in the negative nilpotent part of the enveloping algebra $\U(\gog)$, where $\gog$ is the Lie algebra of $G$, and has remarkable compatibility properties with the simple $G$-modules of highest weight. The basis $\B$ provides good bases of simple $G$-modules. And this provides good bases to study the homogenous coordinates algebras of the varieties.\par 
By a result of Lusztig, we know that $w_0$ acts by a permutation on the elements of the bases of the $G$-modules induced by $\B$. To have explicit results we use a combinatorics of $\B$, given in terms of string parametrization and Lusztig parametrization. These parametrizations depend on a choice of a reduced decomposition of $w_0$. In the case where $G=SL_n(\C)$, for a convenient choice of the decomposition of $w_0$, this combinatorics is the same as the combinatorics given in terms of Young tableaux. In this case, the action of $w_0$ is given by the involution of Sch\"utzenberger described on the tableaux in \cite{schutzenberger1}, and we have explicit formulas. This was done in \cite{berensteinzelevinsky1}. We generalize these results (Corollary \ref{maformuletrop} and Corollary \ref{formuletropeta}) to any group $G$ and to any choice of a reduced decomposition. A part of our results was already announced in \cite{moi}, and applied in \cite{CMMG}.\par 
The generalized Sch\"utzenberger involution is understood via the geometric lifting, i.e a geometric version of the canonical basis which gives a combinatorics of totally positive subvarieties in $G$. We give (Theorem \ref{ThTropZeta}) a geometric analogue of the Sch\"utzenberger involution in the totally positive subvarieties of $G$. The formulas in the geometric version can be easily computed. These formulas are closely related to similar formulas in the algebraic version by a "tropicalization" application. We strongly use the results of \cite{berensteinfominzelevinsky} and \cite{berensteinzelevinsky2}.\par 
This paper is organized as follows. Section~1 provides a construction of the canonical basis and their parametrizations. It also recalls the compatibility property with the simple highest weight $G$-modules. In Section~2, we define the action of $w_0$ on the modules and we give its geometric analogue. We obtain explicit formulas in terms of parametrizations of the canonical basis. In Section~3, we recall the constructions of degenerations of the flag variety and the Schubert varieties due to \cite{Philippe}. We then construct semitoric degenerations of the Richardson varieties.\par 
I wish to thank P.Caldero for introducing me the statement of this problem and for his constant help. 
I am deeply grateful to P.Baumann for enlightening discussions and for his help with Lemma~\ref{lemmepierre}.

\section{Notations and preliminaries}
\subsection{Main data and notations}
Let $G$ be a semisimple simply connected complex Lie group. Fix a torus $T$ and a Borel subgroup $B$
of $G$ such that $T\subset B \subset G$. Let $N$ be the unipotent radical of $B$. Denote by $B^-$ the opposite Borel subgroup and $N^-$ its unipotent radical. The complex Lie algebras associated to $G$, $T$, $N$, $N^-$ will be denoted by $\gog$, $\gh$, $\gn$, $\gn^-$ respectively. 
There is a triangular decomposition $\gog=\gn^-\oplus \gh\oplus \gn $. Let $\lbrace\alpha _i
\rbrace_{1\leq i\leq n}$ be the set of simple roots corresponding to this decomposition, where $n$ is the rank of $\gog$. This set provides a basis of the dual vector space $\gh^*$. The simple coroots in $\gh$ are denoted by $\lbrace \alpha _i ^{\vee}\rbrace_{1\leq i\leq n}$. The weight lattice $P:=\{ \lambda \in \gh^*, \: \lambda(\alpha_i^{\vee}) \in \Z, \forall \:1\leq i\leq n \}$ is generated by the fundamental weights $\varpi _i$, $1\leq i\leq n$, defined such that $\varpi _j(\alpha _i^{\vee })=\delta _{i,j}$. Let $P^+:=\sum_i\N.\varpi _i$ be the semigroup of integral dominant weights. The natural bilinear form on $\gh^* \times \gh$ is denoted by $\langle \; , \; \rangle$. The Cartan matrix associated to 
$\gog$ is $(a_{ij})_{1\leq i,j\leq n}$; one has $a_{ij}=\langle \alpha_j,\alpha_i ^{\vee}\rangle$. Recall that $a_{ii}=2$, $a_{ij}\leq 0$, for all $1\leq i\neq j\leq n$, and there exist nonnegative integers $(d_i)_{1\leq i\leq n}$ such that $d_ia_{ij}=d_ja_{ji}$.

\subsection{Weyl group and reduced words}
The Weyl group $W$ is the subgroup of $\textup{End}(\gh^*)$ generated by the reflexions $s_i$, $1\leq i\leq n$, such that $s_i(\lambda)=\lambda-\langle \lambda, \alpha_i^\vee \rangle \alpha_i$, $\forall \: \lambda \in \gh^*$. We identify $s_i$ with its adjoint so we also have $s_i(h)=h-\langle \alpha_i,h\rangle\alpha_i^{\vee}$, $\forall \: h \in \gh$. The form $\langle \; , \; \rangle$ is $W$-invariant.
 A \textit{reduced word for} $w\in W$ is a finite sequence of indices 
$\ii=(i_1,\cdots, i_l)$ such that $w=s_{i_1}\cdots s_{i_l}$ and the length $\ell (w):=\ell $ is the shortest possible length. Let $w_0$ be the unique element of $W$ with maximal length; set $N:=\ell (w_0)$. Reduced words for $w=w_0$ will be called \textit{reduced words} for short.
The involution $i \mapsto i^* $ of the set $\lbrace 1,\dots,n\rbrace$ is defined by $w_0(\alpha_i)=-\alpha_{i^*}$. Given a reduced word 
$\ii=(i_1,\cdots, i_N) $, we set $\ii^*:=(i_1^*,\cdots, i_N^*)$. It is clear that $\ii^*$ is also a reduced word. Given $\lambda \in P^+$, we set $\lambda ^*:=-w_0(\lambda )$.

\subsection{PBW-bases} \label{quant1}
Now, let us introduce the quantum enveloping algebras. They will be useful in Sections \ref{quant1}, \ref{quant2} and \ref{quant3} for the definition of the canonical basis and there parametrizations. After these sections, we will only consider the classical algebras which are the specializations at $q=1$.
Let $q$ be an indeterminate. The quantum enveloping algebra $\U_q(\gog)$ of $\gog$, over $\C(q)$, is defined with generators $E_i$, $F_i$, $K_i$, $1\leq i\leq n$ and quantum Serre relations. We also have a triangular decomposition $\U_q(\gog)=\U_q(\gn^-)\otimes \U_q(\gh)\otimes \U_q(\gn)$. 
One can construct bases of $\U_q(\gn)$ called Poincar\'e-Birkhoff-Witt type bases as follows.\par 
For all $1\leq i\leq n$ and all $ k \in \N$, we set $q_i:=q^{d_i}$,   $[k]_{i}:=\frac{q_i^{k}-q_i^{-k}}{q_i-q_i^{-1}}$ and:
 $$E_i^{(k)}:=\dfrac{1}{[k]_i[k-1]_i\cdots [1]_i}E_i^k, \; F_i^{(k)}:=\dfrac{1}{[k]_i[k-1]_i\cdots [1]_i}F_i^k$$
 Define automorphisms $T_i$, $1\leq i\leq n$ of $\U_q(\gog)$ as follows:
\begin{eqnarray*}
 T_i(K_{j}) &=&K_jK_i^{-a_{ij}},\; 1\leq j\leq n \\
 T_i(E_i)&=&-K_{i}^{-1}F_i,\hskip 5mm T_i(F_i)=-E_iK_i, \\
 T_i(E_j)&=&\sum_{k+l=-a_{ij}}
(-1)^{k}{q_i}^{-k}E_i^{(k)}E_jE_i^{(l)}, \; 1\leq j\neq i \leq n\\
T_i(F_j)&=&
\sum_{k+l=-a_{ij}}
q_i^{k}F_i^{(l)}F_jF_i^{(k)}, \; 1\leq j\neq i \leq n
\end{eqnarray*}
One can check the compatibility with Serre's relations.
Now fix a reduced word ${\bf i}=(i_1,\ldots,i_N)$. For all
$k$, $1\leq k\leq N$, set $\beta_{\ii,k}:=s_{i_1}\ldots
s_{i_{k-1}}(\alpha_{i_k})$. It is well known that $\{\beta_{\ii,k},\,
1\leq k\leq N\}$ is the set of positive roots and that the ordering 
$$\beta_{\ii,1}<\beta_{\ii,2}<\ldots<\beta_{\ii,N}$$
is a convex ordering on $R^+$.
For all $k$, we define $E^{{\bf i}}_{\beta_{\ii,k}}=T_{i_1}\ldots
T_{i_{k-1}}(E_{i_k})$.
Furthermore, for all $t=(t_1,\cdots t_N)\in\N^N$, we set $E^{{\bf i}}(t):=E_{\beta_{\ii,1}}^{(t_1)}\ldots
E_{\beta_{\ii,N}}^{(t_N)}$, where
$E_{\beta_{\ii,k}}^{(t_k)}:=\frac{1}{[t_k]_{\beta_{\ii,k}}!}E_{\beta_{\ii,k}}^{t_k}$.
The set $\{E^{\ii}(t), \, t\in\N^N\}$ is the so-called Poincar\'e-Birkhoff-Witt type basis of 
$\U_q(\gn)$ associated to the reduced word ${\bf i}$. In the same way, we define a Poincar\'e-Birkhoff-Witt basis $\{F^{\ii}(t), \, t\in\N^N\}$ of $\U_q(\gn^-)$.
\smallskip

\subsection{Canonical/global basis and its Lusztig parametrization} \label{quant2}
Lusztig, and independantly 
Kashiwara \cite{kashiwara1}, constructed a basis called canonical (or global) basis of the nilpotent part 
$\U_q(\gn^-) $ which have good compatibility properties with the $\gog$-modules. Following Lusztig's construction, let us introduce the "bar" automorphism of $\U_q(\gog)$ over $\C$, denoted $ \,\bar{ } \,$ and defined by:
$$\bar E_i=E_i, \,\bar K_i=K_i^{-1},\,\bar F_i=F_i, \,\bar q=q^{-1},\,
1\leq i\leq n$$

\begin{prop} \cite{lusztig2}
Let ${\bf i}$ be a reduced word. For all $t$ $\in \Z^N_{\geq 0}$, there exists a unique element $b=b_\ii(t)$ in $\U_q(\gn^-)$ such that $\bar b=b$ and $b-F^{{\bf i}}(t)\in
q^{-1}\sum \Z[q^{-1}]F^{{\bf i}}(t')$. The set $\B:=\{ b_\ii(t),t
\in\Z^N_{\geq 0}\}$ does not depend on the choice of the reduced word $\ii$. Moreover $\B$ is a basis of $\U_q(\gn^-)$.
\end{prop}
The set $\B$ as above is namely the canonical basis of $\U_q(\gn^-)$.\\
Given a reduced word $\ii $, the map $t\mapsto b=b_\ii(t)$ is a bijection
from $\Z^N_{\geq 0}$ to $\B$, it gives a  parametrisation of the canonical basis that we call \textit{Lusztig's parametrization}. \par

\subsection{Kashiwara operators and string parametrization}\label{quant3}
Kashiwara's operators acting on the canonical basis may be defined as follows. For all $i$ in $\{1, \cdots n\}$, there exists a unique injective map $\tilde f_i$:~$\B\rightarrow\B$, such that if ${\bf i}$ starts with $i_1=i$,  then:
$$ \tilde f_i(b_\ii(t_1,t_2,\ldots,t_N))=b_\ii(t_1+1,t_2,\ldots,t_N).$$
We also define  $\tilde e_i$ : $\B\rightarrow\B\cup\{0\}$ by $\tilde
e_i(b)=b'$ if there exists $b'$ such that $\tilde f_i(b')=b$ and $\tilde
e_i(b)=0$ otherwise.
We set $\varepsilon_i(b)=$Max$\{k\,| \,\tilde{e}^k_i(b)\not=0\}$. \par

The \textit{string parametrization} of an element $b \in \B$ associated to a reduced word
$\ii=(i_1,\cdots, i_N)$ is the $N$-tuple $c_{\ii}(b):=(t_1,t_2\cdots,t_N)$ defined recursively by

$$t_1=\varepsilon_{i_1}(b),\,t_2=\varepsilon_{i_2}(\tilde
e_{i_1}^{t_1}(b)),\ldots,
t_N=\varepsilon_{i_N}(\tilde e_{i_{N-1}}^{t_{N-1}}\ldots \tilde
e_{i_1}^{t_1}(b)).$$
We denote by $\mathcal{C}_\ii$ the image of $\B$ in $\Z^N_{\geq 0}$ under the map $c_{\ii}$.

\begin{prop} \cite{kashiwara1} \label{opeKashi}
Let $\ii=(i_1,\cdots, i_N) $ be a reduced word, and let $b$ be an element of $ \B$
with string parameter $c_{\ii}(b)=(t_1,t_2\cdots,t_N)$. One has
\item[(i)] $\tilde f_{i_1}^{t_1}\ldots \tilde f_{i_N}^{t_N}(1)=b$, 
\item[(ii)] $c_{\ii}(\tilde f_{i_1}(b))=(t_1+1,t_2,\cdots,t_N)$
\end{prop}

\subsection{Transition maps $R_{\pm \ii}^{\pm \ii'}$} 
Let us now introduce the various reparametrization maps. Let $\ii$ and $\ii'$ be reduced words, define: 

$$R_{{\bf i}}^{{\bf i}'}=(b_{{\bf i}'})^{-1}\circ
b_\ii\;:\;\N^N\rightarrow\N^N,$$
$$R_{-{\bf i}}^{-{\bf i}'}=c_{{\bf i}'}\circ(c_\ii)^{-1}\;:\;\mathcal{C}
_{\ii} \rightarrow{\mathcal C}_{{\bf i}'},$$
$$R_{-{\bf i}}^{{\bf i}'}=(b_{{\bf
i}'})^{-1}\circ(c_\ii)^{-1}\;:\;\mathcal{C} _{\ii} \rightarrow\N^N,$$
$$R_{{\bf i}}^{-{\bf i}'}=c_{{\bf i}'}\circ b_\ii
\;:\;\N^N\rightarrow{\mathcal C}_{{\bf i}'}.$$

\begin{ex}
In the case $G=SL_3$, there are exactly two reduced words, namely $\ii=(1,2,1)$ and $\ii'=(2,1,2)$. The map
$R_{{\bf i}}^{{\bf i}'}$ was calculated by Lusztig, see \cite{lusztig1}. If
$b_\ii(a,b,c)=b_{\ii'}(a',b',c')$, then:
 
\begin{equation}\label{calculLusztig}\left\lbrace \begin{array}{lcl}
a'& = & b+c-\min(a,c)\\
b'& = & \min(a,c)\\
c'& = & a+b-\min(a,c) \end{array}
\right.
\end{equation} 

\end{ex}

The methods of computation and explicit formulas of all the previous maps are given in \cite{berensteinzelevinsky2}, we will recall them in section \ref{ThmTropR}.

\subsection{Canonical basis in the modules}
Given a weight $\lambda $ in $P^+$, the Weyl module denoted by $V(\lambda )$
is a simple finite dimensional $\U(\gog)$-module with highest weight $\lambda $. From now on, we fix for any $\lambda \in P^+$, a highest weight vector $v_\lambda $ and a lowest weight vector $v^{low}_\lambda $ in every $V(\lambda )$. One has $V(\lambda )=\U(\gn^-).v_\lambda =\U(\gn).v^{low}_\lambda $. It is known that the module $V(\lambda )$ satisfies the Weyl character formula.
Let $w$ be an element in $W$, fix an extremal vector $v_{w\lambda }$ in $V(\lambda )$ of weight $w\lambda $. We introduce the Demazure module $V_w(\lambda ):=\U(\gn).v_{w\lambda} $ which is a $\U(\gb)$-submodule of $V(\lambda )$.\par
The canonical basis and the above modules are compatible, by \cite{kashiwara1} and \cite{littelmann}.

\begin{thm} \label{ThmBaseModule}
One has:
\begin{enumerate}
\item If $\B(\lambda) :=\{b\in \B, \, bv_\lambda \not=0\}$, then 
$\B(\lambda) v_\lambda $ is a basis of $V(\lambda )$.
\item There exists a subset $\B_w$ of $\B$, which does not depend on $\lambda
$, such that $\B_wv_\lambda $ generates $V_w(\lambda )$.
\end{enumerate}
\end{thm}
We will use abreviation $b$ instead of $bv_\lambda$ and $\B(\lambda)$ instead of $\B(\lambda)v_\lambda$ when no confusion occurs. Denote $\B_w(\lambda):=\B(\lambda)\cap \B_w$.\par 
Now, we may suppose that $v^{low}_\lambda $ and $v_{w\lambda }$ belong to $\B(\lambda)$.\par
We still denote by $\tilde{e}_i$ and $\tilde{f}_i$ the Kashiwara operators define from $\B(\lambda)$ to $\B(\lambda)\cup \{0\}$ by $\tilde{e}_i(bv_\lambda)=\tilde{e}_i(b)v_\lambda$ and 
$\tilde{f}_i(bv_\lambda)=\tilde{f}_i(b)v_\lambda$.

\subsection{Examples in the $\texttt{A}_n$ case}
In this section, we study the case where $G=SL_{n+1}$. In this case the Weyl group is isomorphic to the group of 
permutations $\mathfrak{S}_n$. The element $w_0$ has length $n(n+1)/2$ and the special reduced word  $\ii=(1,2,1,\cdots, n, n-1, \cdots, 2,1)$ will be called \textit{standard reduced word}. Recall that in the case $G=SL_{n+1}$, one has a combinatoric model of Young tableaux. The Lusztig parametrization and the string parametrization generalize this combinatory. They coincide when the parametrizations are considered with the standard reduced word.

\begin{defn} 
Let $\lambda =\lambda _1\varpi _1+ \lambda _2\varpi _2+\cdots \lambda _n\varpi _n$ be an element of $P^+$. The Young tableau of shape $\lambda $ is a collection of boxes, arranged from left to right, from $\lambda _n$ columns with $n$ boxes to $\lambda _1$ columns with one boxes. Further, a tableau filled with entries in $\{1,\cdots, n+1\}$ such that the entries increase across each row and strictly increase down each column, is called a semi-standard Young tableau of shape ~$\lambda$.

\end{defn}

%\begin{ex}$ $\\

%\noindent
%\upshape
%Soit $n\geq 3$, \\

%$\yng(4,3,1)$ est un tableau de Young de forme $\varpi _1+2\varpi _2+\varpi _3$.
%\end{ex}

\begin{ex}
Let $n\geq 3$, \\

$\young(1223,234,4)$ is a semi-standard Young tableau of shape $\varpi _1+2\varpi _2+\varpi _3$
\end{ex}

Denote by $Y(\lambda)$ the set of all semi-standard Young tableaux of shape $\lambda$.
One knows (see \cite{kashiwara1}) that the set $Y(\lambda)$ gives a parametrization of the canonical basis $\B(\lambda)$. Let us precise this fact.\par
Let $T$ be a tableau $Y(\lambda)$ and denote by $b_T$ the element associated to the tableau $T$. Let 
$\ii$ be the standard reduced word, introduce Lusztig's parameters of $b_T$:
$$(t_{11}, t_{12}, t_{22}, t_{13}, t_{23}, t_{33} \cdots, t_{1n}, \cdots, t_{nn}):=b_\ii^{-1}(b_T)$$
and the string parameters of $b_T$:  
$$(c_{11}, c_{22}, c_{12}, c_{33}, c_{23}, c_{13} \cdots, c_{nn}, \cdots, c_{1n}):=c_\ii(b_T)$$

The link between all theses parametrizations is the following.

\begin{prop} \cite{berensteinzelevinsky1}
One has,
\begin{equation} \begin{array}{lcl}
t_{ij}&=& (\text{the number of $j+1$ in the $i$-th row of T}), \;  1\leq i \leq j \leq n\\
c_{ij}&=& (\text{the number of $j+1$ in the $i$ first rows of T}), \; 1\leq i \leq j \leq n \end{array}
\end{equation}
\end{prop}

One can deduce explicit formulas for the maps $R_\ii^{-\ii}$ and $R_{-\ii}^\ii$ in this special case.

\begin{cor} 
One has:
\begin{equation} \label{chgtlineaire} \begin{array}{lcl}
t_{ij}&=&c_{i+1,j}-c_{ij}, \; 1\leq i \leq j \leq n\\
c_{ij}&=&t_{1j}+t_{2j}+\cdots +t_{ij},  \; 1\leq i \leq j \leq n \end{array}
\end{equation}
\end{cor}

\begin{ex} Figure \ref{Gloup2} describes the canonical basis $\B(\varpi_1+\varpi_2)$ in the case $A_2$. For each element of the basis, we give the Young tableau associated, the string parameters and the Lusztig's parameters for the standard reduced word $\ii=(1,2,1)$. In this figure, the simple arrows between two elements $b\longrightarrow b'$ mean that $b'=\tilde{f}_1(b)$, and the double arrows $b\Longrightarrow b'$ mean that $b'=\tilde{f}_2(b)$.

%\begin{figure}[!ht]
%\begin{pspicture}(13,14)(0,-1)
% \put(4,-0.25){\shortstack{$(1,2,1)$\\ $(1,1,1)$}}
% \put(5.5,0){$\young(23,3)$}
 
% \psline{>}(3.9,0.5)(1.5,2.5)

 % \put(0.3,2.7){\shortstack{$(0,2,1)$\\ $(0,1,1)$}}
 %\put(1.8,2.95){$\young(13,3)$}
 
 %\psline[doubleline=true]{>}(1.5,3.5)(1.5,5.9)

 % \put(0.3,6.1){\shortstack{$(0,1,1)$\\ $(0,1,0)$}}
% \put(1.8,6.35){$\young(13,2)$}
 
% \psline[doubleline=true]{>}(1.5,6.9)(1.5,9.3)

 %\put(0.3,9.5){\shortstack{$(1,0,0)$\\ $(1,0,0)$}}
%\put(1.8,9.75){$\young(12,2)$}
 
 %\psline{>}(1.55,10.35)(3.9,12.3)
 %\psline[linewidth=.01]{<-}(2.225,11.325)(2.725,11.74)
 %\put(2.25,11.825){$\tilde{f}_1$}
 
 %\put(3.7,12.45){\shortstack{$c_\ii(0,0,0)$\\ $\hspace{-0.24cm} b^{-1}_\ii(0,0,0)$}}
 %\put(5.5,12.8){$\young(11,2)$}
 
 %\psline[doubleline=true]{->}(6.4,12.3)(8.8,10.35)
 %\psline[linewidth=.01]{<-}(8.1,11.325)(7.6,11.74)
 %\put(7.65,11.825){$\tilde{f}_2$}
 
 %\put(7.6,9.5){\shortstack{$(0,1,0)$\\ $(0,0,1)$}}
 %\put(9.1,9.75){$\young(11,3)$}
  
 %\psline{>}(8.8,6.9)(8.8,9.3)
 
 %\put(7.6,6.1){\shortstack{$(1,1,0)$\\ $(1,0,1)$}}
 %\put(9.1,6.35){$\young(12,3)$}
 
 %\psline{>}(8.8,3.5)(8.8,5.9)
 
% \put(7.6,2.7){\shortstack{$(2,1,0)$\\ $(2,0,1)$}}
 %\put(9.1,2.95){$\young(22,3)$}
 
 %\psline[doubleline=true]{->}(8.8,2.5)(6.4,0.6)

%\end{pspicture}
%\caption{Parametrizations of the basis $\B(\varpi_1+\varpi_2)$}
%\label{Gloup}
%\end{figure}

\end{ex}

\section{Action of $w_0$ and geometric lifting}
\subsection{Modules twisted by automorphism} \label{definvol}
Let us consider the three automorphisms of $\U(\gog)$ defined on the generators by:
$$\phi (E_i)=F_i, \hspace{0.5cm} \phi (F_i)=E_i, \hspace{0.5cm} \phi
(H_i)=-H_i$$
$$\delta (E_i)=E_{i^*}, \hspace{0.5cm} \delta (F_i)=F_{i^*},
\hspace{0.5cm} \delta (H_i)=H_{i^*}$$
$$\eta (E_i)=F_{i^*}, \hspace{0.5cm} \eta (F_i)=E_{i^*}, \hspace{0.5cm}
\eta (H_i)=-H_{i^*}$$
Notice that the automorphism $\eta$ coincides with the action of $w_0$ up to a multiplicative constant.\par
In the sequel, let us fix a dominant weight $\lambda \in P^+$.\par
Given an automorphism $\chi $ of $\U(\gog)$, one can define the twisted module $V(\lambda )^\chi $ as the vector space $V(\lambda )$ with the following action: $u*v=\chi (u)v$, $u\in
\U(\gog)$, $v\in V(\lambda )$. The module $V(\lambda )^\chi$ is simple since $V(\lambda )$ is simple. And one has  $V(\lambda )^\chi \simeq V(\lambda ^\chi )$, for a certain $\lambda ^\chi \in P^+$. The automorphism $\chi $ leads an isomorphism of vector spaces 
$\chi _\lambda :V(\lambda ) \longrightarrow V(\lambda ^\chi )$,
verifying $\chi _\lambda (uv)=\chi (u)\chi _\lambda (v)$. Such isomorphism  is unique up to  multiplicative constant by Schur's lemma. \par 

Let us describe these isomorphisms in the cases where $\chi=\eta, \, \delta$ and $\phi$. The isomorphism $\eta _\lambda :V(\lambda )\rightarrow V(\lambda^{\eta} )$ satisfies $\eta _\lambda (u.v)=\eta (u)\eta _\lambda (v)$, thus one has:
$$F_i\eta _\lambda (v_{\lambda})=\eta(E_{i^*})\eta _\lambda (v_{\lambda})=\eta _\lambda (E_{i^*}v_{\lambda})=0, \; \forall 1\leq i \leq n $$
The vector $\eta_{\lambda}(v_{\lambda})$ is therefore a lowest weight vector in the corresponding twisted module $V(\lambda^{\eta})$. Let us determine the weight of $\eta_{\lambda}(v_{\lambda})$. For all $1\leq i \leq n$, one has:
$$ H_i\eta_{\lambda}(v_{\lambda})=\eta_\lambda(-H_{i^*}v_{\lambda})=\eta_\lambda(-\langle \lambda, \alpha_{i^*}^{\vee}\rangle v_{\lambda})=-\langle \lambda, \alpha_{i^*}^{\vee}\rangle \eta_{\lambda}(v_{\lambda})=\langle w_0(\lambda), \alpha_{i}^{\vee}\rangle \eta_{\lambda}(v_{\lambda})$$
Hence, $\eta_{\lambda}(v_{\lambda})$ is a lowest weight vector of weight $w_0(\lambda)$. One deduces that $V(\lambda^{\eta})\simeq V(\lambda)$ and that $\eta_{\lambda}(v_{\lambda})$ is proportional to $v^{low}_\lambda$.\par 
From now on, set  $\eta_{\lambda}(v_{\lambda})=v^{low}_\lambda$. To summarize, one has:
$$\eta _\lambda :V(\lambda )\rightarrow V(\lambda ) \text{ with }
\eta_\lambda (uv_\lambda )=\eta(u)v^{low}_\lambda, \; \forall u\in \U(\gog)$$

In the same way, the automorphisms $\phi$ and $\delta$ induce the following isomorphisms of vector spaces (normalized by the choice of the image of $v_\lambda$):

$$\phi _\lambda :V(\lambda )\rightarrow V(\lambda ^*) \text{ with } \phi
_\lambda (uv_\lambda )=\phi(u)v^{low}_{\lambda ^* }, \; \forall u\in \U(\gog)$$
$$\delta _\lambda :V(\lambda )\rightarrow V(\lambda ^*) \text{ with }
\delta _\lambda (uv_\lambda )=\delta(u)v_{\lambda ^* }, \; \forall u\in \U(\gog)$$

The isomorphism $\phi_{\lambda}$ is compatible with the canonical basis in the following sense:

\begin{prop}\cite[\S 21]{lusztig1} \label{baseca_et_autom} One has:
\item[(i)] $ \hspace{0.2cm} \phi _\lambda (\B(\lambda )v_ \lambda )=\B(\lambda ^* )v_{\lambda ^* }$
\item[(ii)] $  \; \forall \,1\leq i\leq n, \;, \forall b \in \B(\lambda), \; \tilde{e}_i\, \phi _\lambda(b)=\phi _\lambda \tilde{f}_i(b)\, $
\end{prop}

It is clear from the definitions that $\delta (b_\ii(t))=\delta (b_{\ii^*}(t))$, thus we also have 
$\delta_\lambda (\B(\lambda )v_ \lambda )=\B(\lambda ^* )v_{\lambda ^* }$. It is also clear that 
$\eta_\lambda=\phi_{\lambda} \delta_{\lambda}$, and thus $\eta_\lambda (\B(\lambda )v_ \lambda )=\B(\lambda )v_{\lambda }$.

\subsection{Sch\"utzenberger involution}
Fix a dominant weight $\lambda \in P^+$.
The isomorphism $\eta _\lambda $ generalizes the Sch\"utzenberger involution defined in the case $G=SL(n)$ in terms of Young tableaux. Sch\"utzenberger described (see \cite{schutzenberger1}) an involution $S:Y(\lambda)\rightarrow Y(\lambda)$ with an algorithm called "jeu de taquin". \\
In the case where $G=SL(n)$, the link between 
$\eta_\lambda$ and $S$ is the following:

\begin{prop} \cite{berensteinzelevinsky1} Given $T\in Y(\lambda)$ and the associated element $b_T \in \B(\lambda)$, one has:
$$\eta_\lambda (b_T)=b_{S(T)}$$

\end{prop}

\begin{ex}
Figure \ref{Gloup2} describes the Sch\"utzenberger involution on the basis of the $\mathfrak{sl}_2$-module $V(\varpi_1+\varpi_2)$.

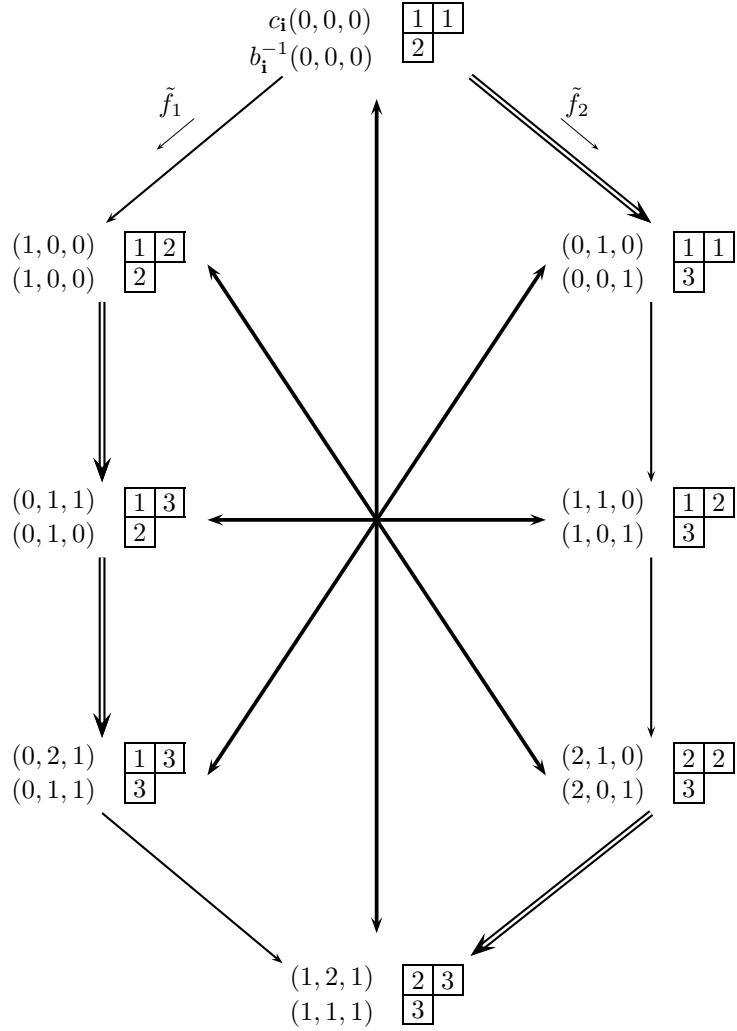
\begin{figure}[!ht]
\begin{pspicture}(13,14)(0,-1)
 \put(4,-0.25){\shortstack{$(1,2,1)$\\ $(1,1,1)$}}
 \put(5.5,0){$\young(23,3)$}
 
 \psline{>}(3.9,0.5)(1.5,2.5)

  \put(0.3,2.7){\shortstack{$(0,2,1)$\\ $(0,1,1)$}}
 \put(1.8,2.95){$\young(13,3)$}
 
 \psline[doubleline=true]{>}(1.5,3.5)(1.5,5.9)
 %\psline(1.6,3.5)(1.6,5.9)
 
  \put(0.3,6.1){\shortstack{$(0,1,1)$\\ $(0,1,0)$}}
 \put(1.8,6.35){$\young(13,2)$}
 
 \psline[doubleline=true]{>}(1.5,6.9)(1.5,9.3)
 %\psline(1.6,6.9)(1.6,9.3)
 
  \put(0.3,9.5){\shortstack{$(1,0,0)$\\ $(1,0,0)$}}
 \put(1.8,9.75){$\young(12,2)$}
 
 \psline{<-}(1.55,10.35)(3.9,12.3)
 \psline[linewidth=.01]{<-}(2.225,11.325)(2.725,11.74)
 \put(2.25,11.825){$\tilde{f}_1$}
 
 \put(3.7,12.45){\shortstack{$c_\ii(0,0,0)$\\ $\hspace{-0.24cm} b^{-1}_\ii(0,0,0)$}}
 \put(5.5,12.8){$\young(11,2)$}
 
 \psline[doubleline=true]{->}(6.4,12.3)(8.8,10.35)
 \psline[linewidth=.01]{<-}(8.1,11.325)(7.6,11.74)
 \put(7.65,11.825){$\tilde{f}_2$}
 
 \put(7.6,9.5){\shortstack{$(0,1,0)$\\ $(0,0,1)$}}
 \put(9.1,9.75){$\young(11,3)$}
  
 \psline{>}(8.8,6.9)(8.8,9.3)
 
 \put(7.6,6.1){\shortstack{$(1,1,0)$\\ $(1,0,1)$}}
 \put(9.1,6.35){$\young(12,3)$}
 
 \psline{>}(8.8,3.5)(8.8,5.9)
 
 \put(7.6,2.7){\shortstack{$(2,1,0)$\\ $(2,0,1)$}}
 \put(9.1,2.95){$\young(22,3)$}
 
 \psline[doubleline=true]{->}(8.8,2.5)(6.4,0.6)

 %fleches de l'action de w_0:
 
\psline[linewidth=.05]{<->}(5.15,12)(5.15,0.9)
\psline[linewidth=.05]{<->}(2.9,9.8)(7.4,3)
\psline[linewidth=.05]{<->}(2.9,3)(7.4,9.8) 
\psline[linewidth=.05]{<->}(2.9,6.4)(7.4,6.4)
 
\end{pspicture}
\caption{Canonical basis $\B(\varpi_1+\varpi_2)$ and Sch\"utzenberger involution.}
\label{Gloup2}
\end{figure}

\end{ex}

\subsection{Geometric lifting and tropicalization}
In the sequel, we wish to compute explicit formulas for the isomorphism $\eta_\lambda$ in terms of parametrizations of the canonical basis. We first consider the isomorphism $\phi _\lambda $, and then we use the composition $\eta_\lambda =\delta _{\lambda ^* }\phi _\lambda $.\par
Our goal is to give explicit formulas for the application
$b_{\ii}^{-1}\phi _\lambda c_{\ii}^{-1}$ which expresses Lusztig's parameters $t'=(t'_1,\cdots, t'_N)$ of
the element $\phi _\lambda (b) \in \B(\lambda^*)$ in terms of the string parameters $t=(t_1,\cdots ,t_N)$ of an element $b \in \B(\lambda)$.\par 
For this end we use the methods and the results of \cite{berensteinfominzelevinsky} and \cite{berensteinzelevinsky2} on geometric lifting and tropicalization.\par 
We consider semifield structures, i.e commutative multiplicative groups equipped with an additive law which is associative, commutative and distributive on the product.\par 
The main example of semifield is the set of integers $\Z$ endowed with the operations: $a\oplus b:=\min (a,b)$, $a\odot b:=a+b$, $a, b \in \Z$, it is called \textit{tropical structure} of $\Z$.  
It induces a semifield structure on the set $\mathcal{F}(\Z^N_{\geq 0}, \Z)$ of all maps from
$\Z^N_{\geq 0}$ to $ \Z$ with the operations:
$$f\odot g: (t_1, \cdots, t_N) \mapsto f (t_1, \cdots, t_N)+g(t_1,
\cdots, t_N), $$
$$ f\oplus g:(t_1, \cdots, t_N)\mapsto \min (f(t_1, \cdots, t_N), g(t_1, \cdots, t_N))$$
for all $f,g\in \mathcal{F}(\Z^N_{\geq 0}, \Z)$.\par 
We denote by $p_i \in \mathcal{F}(\Z^N_{\geq 0}, \Z)$ $1\leq i\leq N$, the projections defined by 
 $p_i(t_1, \cdots, t_N) =t_i$. Let $\Q_{>0}(t_1,\cdots,t_N)$ be the set of \textit{ rational
  subtraction-free expressions} in inderterminates $t_1,\cdots, t_N$.
  This set is a semifield (the smallest one containing the indeterminates $t_1,\cdots, t_N$) for the
 usual laws + and $\times $. The \textit{tropicalization} is defined by:

\begin{thm}[\cite{berensteinfominzelevinsky}]
There exists a unique homomorphism of semifields, denoted by $ [.]_\tr$, such that:
\begin{eqnarray*}
 [.]_\tr :\Q_{>0}(t_1,\cdots,t_N) & \rightarrow &\mathcal{F}(\Z^N_{\geq 
0},\Z), \\
 t_i & \mapsto & p_i, \; 1\leq i\leq n,
\end{eqnarray*}

\end{thm}

In other words, if $f(t_1,\cdots,t_N)$ is a subtraction-free expression, the tropicalization $[f]_\tr(t_1,\cdots,t_N)$
is the expression obtained from $f(t_1,\cdots,t_N)$ by changing the  + in $\min$, the $\times $ 
in +, et the $\div $ in - . Let us give an example taken from \cite{berensteinfominzelevinsky}.

\begin{ex}
Consider the expression $f(t_1,t_2):=t_1^2-t_1t_2+t_2^2$. 
It is a rational subtraction-free expression since 
$f(t_1,t_2)=\frac{t_1^3+t_2^3}{t_1+t_2}$. One has 
$[f]_\tr(t_1,t_2)=\min(3t_1,3t_2)-\min(t_1,t_2)=\min(2t_1,2t_2)$.
\end{ex}
The element $f$ is called \textit{geometric lifting} of $[f]_\tr$.\par
Notice that geometric lifting is not uniquely determined by $[f]_\tr$ as we can see in the above example.
Indeed, 
$\frac{t_1^3+t_2^3}{t_1+t_2}$ and $t_1^2+t_2^2$ are both geometric lifting of $\min(2t_1,2t_2)$.\par

To finish this subsection, let us introduce the following notation.\\
Let $f_1, f_2,\cdots, f_n \in \Q_{>0}(t_1,\cdots,t_N)$, we set
$$[(f_1, f_2,\cdots, f_n)]_\tr:=([f_1]_\tr,[f_2]_\tr,\cdots,[f_n]_\tr)$$

\subsection{Totally positive subvariety of $G$}
For any $1\leq i\leq n$, denote by $\varphi _i:SL_2\hookrightarrow G$
the natural injective map corresponding to the simple root $\alpha _i$.
Consider the one-parameter subgroups of $G$ defined by
$$x_i(t)=\varphi _i\left(
\begin{array}{ccc}
1&t\\
0&1
\end{array}
\right),
\hspace{0.5cm} y_i(t)=\varphi _i\left(
\begin{array}{ccc}
1&0\\
t&1
\end{array}
\right)
,\hspace{0.5cm}t\in \C$$
$$t^{\alpha _i^\vee}=\varphi _i\left(
\begin{array}{ccc}
t&0\\
0&t^{-1}
\end{array}
\right)
,\hspace{0.5cm}t\in \C^*.$$
Clearly, $x_i(t)$'s, (resp. $y_i(t)$, $t^{\alpha _i ^\vee}$) generate $N$, (resp. $N^-$, $T$).
One has the following relations of commutation:
\begin{eqnarray} \label{relcommut}
t^{\alpha _i^\vee}x_j(t')=x_j(t^{a_{ij}}t')t^{\alpha _i^\vee},
\hspace{0.3cm}t^{\alpha _i^\vee }y_j(t')=y_j(t^{-a_{ij}}t')t^{\alpha
_i^\vee}
\end{eqnarray}
One defines two involutive antiautomorphisms of $G$: $x\mapsto x^T$,
called transposition, and $x\mapsto x^\iota $, called inversion, as follows:
$$
\begin{array}{lll}
x_i(t)^T=y_i(t), & y_i(t)^T=x_i(t), & (t^{\alpha _i^\vee})^T=t^{\alpha
_i^\vee}\\[4pt]
x_i(t)^ \iota =x_i(t), & y_i(t)^\iota =y_i(t) ,& (t^{\alpha
_i^\vee})^\iota =t^{-\alpha _i^\vee}
\end{array}
$$

Let $G_0:=N^-TN$ be the set of all elements in $G$ which admit a Gaussian decomposition. Given  $x\in G_0$, the Gaussian decomposition is unique and we will write $x=[x]_-[x]_0[x]_+$, where $[x]_-\in N^-, \: [x]_0 \in T, \; [x]_+ \in N$.\\
Let $G_{\geq 0}$ be the submonoid of $G$ generated by  
$x_i(t)$, $y_i(t)$, $t^{\alpha _i^\vee}$ for all $t>0$.\\
Given a word $\ii=(i_1,\cdots, i_m) $ and a $m$-tuple
$t=(t_1,\cdots, t_m) $ in $\C^m_{\neq 0}$, we set:
\vspace{-0.2cm}
$$x_{\ii}(t):=x_{i_1}(t_1)\cdots x_{i_m}(t_m), \hspace{0.3cm} 
\text{et}
\hspace{0.3cm}
x_{-\ii}(t):=y_{i_1}(t_1)t_1^{-\alpha _{i_1} ^{\vee}} \cdots
y_{i_m}(t_m)t_m^{-\alpha _{i_m} ^\vee} $$

Consider the following reduced double Bruhat cells.
 
\begin{eqnarray}\label{def2cellules}
L^{e,w_0}:=N\cap B_-w_0B_- \hspace{0.5cm} \text{ and } \hspace{0.5cm}
L^{w_0,e}:=Nw_0N\cap B_-
\end{eqnarray}

Denote by $L^{e,w_0}_{>0}$, resp. $L^{w_0,e}_{>0}$, their intersections with $G_{>0}$.
The maps $x_\ii$ and $x_{-\ii}$ parametrize these subvarieties of $G$. More precisely,
\begin{thm}\cite{berensteinzelevinsky2}\label{thm.param.cells}
For any reduced word $\ii$, the map $x_{\ii}$, resp.
$x_{-\ii}$, is a birational isomorphism between $\C_{\not = 0}^N$ and
$L^{e,w_0}$, resp. $L^{w_0,e}$. It restricts to a bijection between
$\R^N_{>0} $ and $L^{e,w_0}_{>0}$, resp. $L^{w_0,e}_{>0}$.
\end{thm}

\subsection{Geometric lifting of the maps $R_{\pm \ii}^{\pm \ii'}$}
In the sequel we will use the notation $(.)^\vee$ that means we consider the analogous maps in the Langlands dual $G^\vee$ of $G$. Recall that the Langlands dual $G^\vee$ is the semisimple Lie group with transposed Cartan's matrix. The simple roots of $G^\vee$ can be naturally identified with the simple coroots of $G$ and conversely. Thus the Weyl groups are naturally identified with each other.\par 
Let us introduce the reparametrization maps
$\tilde {R}_\ii^{\ii'}:=x_{\ii'}^{-1}\circ x_\ii$ and $\tilde
{R}_{-\ii}^{-\ii'}:=x_{-\ii'}^{-1}\circ x_{-\ii}$ from $\C_{\not = 0}^N$ to itself.
An important result from \cite{berensteinzelevinsky2} is that these maps are geometric liftings of
$R_\ii^{\ii'}$ and $R_{-\ii}^{-\ii'}$. More precisely,
\begin{thm} \label{ThmTropR}
The components of $(\tilde {R}_\ii^{\ii'})^\vee (t_1,\cdots ,t_N)$ and $(\tilde
{R}_{-\ii}^{-\ii'})^\vee (t_1,\cdots ,t_N)$ are expressed as rational subtraction-free expressions in $t_1,\cdots ,t_N$, and one has:
$$ (i)\hspace{0.3cm} [(\tilde {R}_\ii^{\ii'})^\vee
]_{\tr}(t)=R_\ii^{\ii'}(t) \hspace{0.5cm} (ii) \hspace{0.3cm} [(\tilde
{R}_{-\ii}^{-\ii'})^\vee ]_{\tr}(t)=R_{-\ii}^{-\ii'}(t) $$
\end{thm}

\begin{ex}Explicit formulas in the case $G=SL_3(\C)$ are:\\
$$ x_1(t)=\left( 
\begin{array}{ccc}
1&t&0\\
0&1&0\\
0&0&1
\end{array}
\right) ,
\; x_2(t)=\left( 
\begin{array}{ccc}
1&0&0\\
0&1&t\\
0&0&1
\end{array}
\right) ,\; t\in \C$$
\vspace*{0.3cm}
$$
t^{\alpha^{\vee}_1}=\left( 
\begin{array}{ccc}
t&0&0\\
0&t^{-1}&0\\
0&0&1
\end{array}
\right) ,
\;t^{\alpha^{\vee}_2}=\left( 
\begin{array}{ccc}
1&0&0\\
0&t&0\\
0&0&t^{-1}
\end{array}
\right),\; t\in \C_{\neq 0}$$
\vspace*{0.3cm}

\noindent  
For any $(t_1,t_2,t_3)\in \C_{\neq 0}^3$ and $(t'_1,t'_2,t'_3)\in \C_{\neq 0}^3$:
\vspace*{0.3cm}
$$x_{121}(t_1,t_2,t_3)=x_1(t_1)x_2(t_2)x_1(t_3)=\left( 
\begin{array}{ccc}
1&t_1+t_3&t_1t_2\\
0&1&t_2\\
0&0&1
\end{array}
\right)$$
\vspace*{0.3cm}
$$
x_{212}(t'_1,t'_2,t'_3)=x_2(t'_1)x_1(t'_2)x_2(t'_3)=\left( 
\begin{array}{ccc}
1&t'_2&t'_2t'_3\\
0&1&t'_1+t'_3\\
0&0&1
\end{array}
\right)$$
\vspace*{0.3cm}

\noindent
If $x_1(t_1)x_2(t_2)x_1(t_3)=x_2(t'_1)x_1(t'_2)x_2(t'_3)$ then:
 
 $$(t'_1,t'_2,t'_3)=(\,\dfrac{t_2t_3}{t_1+t_3},\; t_1+t_3,\; \dfrac{t_1t_2}{t_1+t_3}).$$
Hence,
  $$[(t'_1,t'_2,t'_3)]_{\tr}=(t_2+t_3-\min(t_1,t_3),\: \min(t_1,t_3),\:t_1+t_2-\min(t_1,t_3))$$
 This is precisely the formulas \eqref{calculLusztig} which give the reparametrization $R_{121}^{212}$.
\end{ex}

\begin{ex} 
As in the previous example, we can also compute:
$$x_{-121}(t_1,t_2,t_3)=x_{-212}(t'_1,t'_2,t'_3)$$
One obtains:
$$(t'_1,t'_2,t'_3)=(\,\dfrac{t_2t_3}{t_2+t_1t_3},\; t_1t_3,\; \dfrac{t_2+t_1t_3}{t_3})$$
Hence:
$$[(t'_1,t'_2,t'_3)]_{\tr}=(\:t_2+t_3-\min(t_2,t_1+t_3),\;t_1+t_3, \;\min(t_2,t_1+t_3)-t_3)$$
These formulas give the reparametrization $R_{-121}^{-212}$.
\end{ex}

We can also give a geometric lifting of the maps $R_{-\ii}^{\ii'}$. This geometric lifting is an isomorphism between the subvarieties $L^{w_0,e}$ and $L^{e,w_0}$. To define this isomorphism, we need to introduce a representative of $w_0$ in $G$.\par 
Recall that $W\simeq \textup{Norm}(T)/T$. Fix a representative $\overline{w_0} \in \textup{Norm}(T)$ of $w_0$. In the sequel, the results will not depend on the choice of this representative. \par 
For instance, choose $\overline{w_0}:=\overline{s_{i_1}}\,\overline{s_{i_2}}\cdots\overline{s_{i_N}}$, where $\ii=(i_1,\cdots, i_N)$ is a reduced word and  
$$\overline{s_i}:=\varphi_i\left(
\begin{array}{ccc}
0&-1\\
1&0
\end{array}
\right)=x_i(-1)y_i(1)x_i(-1), \hspace{0.5cm} 1\leq i \leq n$$
We know that $\overline{w_0}$ does not depend on the choice of the reduced word. 
An easy computation shows that 
$\overline{s_i}^T=\overline{s_i}^{-1}$ and $\overline{s_i}^{\iota}=\overline{s_i}$. Furthermore,
\begin{eqnarray}
\overline{w_0}^T=\overline{w_0}^{-1} = \overline{w_0}^\iota= \overline{w_0^{-1}}
\end{eqnarray}
 \noindent
Following \cite{berensteinzelevinsky2}, we define for $x\in G$ 
$$\eta^{w_0,e}(x):=[(\overline{w_0}x^T)^{-1}]_+,$$
and
$$\eta^{e,w_0}(x):=([\overline{w_0}^{-1}x^T]_0[\overline{w_0}^{^{-1}}x^T]_-)^{-1}$$

\begin{thm}{\cite{berensteinzelevinsky2}} \label{ThmTropR2}
\begin{enumerate}
\item The map $\eta^{w_0,e}$ is a birational isomorphism between $L^{w_0,e}$ and $L^{e,w_0}$, which restricts to a bijection from $L^{w_0,e}_{>0}$ to $ L^{e,w_0}_{>0} $; the inverse map is $\eta^{e,w_0}$.
\item The components of $(x_{\ii'}^{-1}\circ \eta^{w_0,e}\circ x_{-\ii})^{\vee}$ are rational subtraction-free expressions, and 
$$R_{-\ii}^{\ii'}(t)=[(x_{\ii'}^{-1}\circ \eta^{w_0,e}\circ x_{-\ii})^{\vee}]_\tr (t)$$
\end{enumerate}
\end{thm}

\begin{ex}
In the case $G=SL_3(\C)$, one has the following explicit formulas:
$$\overline{w_0}=\left( 
\begin{array}{ccc}
0&0&1\\
0&-1&0\\
1&0&0
\end{array}
\right) , \;
x_{-121}(t_1,t_2,t_3)=\left( 
\begin{array}{ccc}
t_1^{-1}t_3^{-1}&0&0\\
t_3^{-1}+t_1t_2^{-1}&t_1t_2^{-1}t_3&0\\
1&t_3&t_2
\end{array}
\right) $$
\vspace*{0.3cm}
$$(\overline{w_0}x_{-121}(t_1,t_2,t_3)^T)^{-1}=\left( 
\begin{array}{ccc}
1& (t_2+t_1t_3)t_3^{-1}&t_1t_3\\
-t_1^{-1}&-t_2t_1^{-1}t_3^{-1}&0\\
t_2^{-1}&0&0
\end{array}
\right)$$

$$=y_1(-t_2t_1^{-1}s^{-1})y_2(-st_2^{-1}t_3^{-1})y_1(-t_3s^{-1})x_1(t_1)x_2(t_3)x_1(t_2t_3^{-1})$$

\noindent
where $s=t_2+t_1t_3$.\\

\noindent
Hence, if $x_{121}(t'_1,t'_2,t'_3)=\eta^{w_0,e}(x_{-121}(t_1,t_2,t_3))$ then
$$(t'_1,t'_2,t'_3)=(t_1,\, t_3, \,t_2t_3^{-1})$$
Hence,
  $$[(t'_1,t'_2,t'_3)]_{\tr}=(t_1,\, t_3, \,t_2-t_3)$$
 This is precisely the formulas \eqref{chgtlineaire} which give the changing of parametrization $R_{-121}^{121}$.
\end{ex}

\subsection{Geometric lifting of $\phi_\lambda$}
Now we fix a dominant weight $\lambda= \lambda_1\varpi_1+\cdots+\lambda_n\varpi_n$ to the end of the section.\\
Given $x \in G$, we set $\zeta (x):=[x^{\iota T}]_{+}$.
Our first observation is 
\begin{prop} \label{FormZeta}
\item[(i)] Let $\ii=(i_1,\cdots, i_N) $ be a reduced word. If $x=x_{-\ii}(t_1,\cdots ,t_N)$ then $\zeta(x)=x_\ii (t_1',\cdots ,t_N')$, where 
\begin{equation}\label{eqFormZeta}
t'_k=t_k^{-1}\prod _{j>k}t_j^{-a_{i_ji_k}}
\end{equation}
\item[(ii)] The application $\zeta$ defines a bijection from $L^{w_0,e}_{>0}$ to $ L^{e,w_0}_{>0} $.
\end{prop}

\noindent
\textbf{Proof}: (i) follows from the relations \eqref{relcommut} and (ii) follows from Theorem~\ref{thm.param.cells} and \eqref{eqFormZeta}.\\

Now we can give a geometric lifting and an explicit formula for $\phi _\lambda $:
\begin{thm}\label{ThTropZeta} 
Let $\ii$ and $\ii'$ be reduced words.
\item[(i)] Components of $(x_\ii^{-1}\circ \zeta \circ x_{-\ii'})^\vee$ are rational subtraction-free expressions.
\item[(ii)] One has,
$$ b_{\ii}^{-1}\phi _\lambda c_{\ii'}^{-1}(t)= [(x_\ii^{-1}\circ \zeta
\circ x_{-\ii'})^\vee ]_\tr(t)+b_{\ii}^{-1}\phi _\lambda (v_\lambda )$$
\end{thm}

As a consequence of the above proposition and the above theorem, we have:

\begin{cor} 
If $(t_1',\cdots ,t_N')=b_{\ii}^{-1}\phi _\lambda c_{\ii}^{-1}(t_1,
\cdots,t_N)$, then
\begin{equation}\label{maformuletrop}
t'_k=l_k-t_k-\sum _{j>k}a_{i_ki_j}t_j
\end{equation}
where $(l_1,\cdots ,l_N):=b_{\ii}^{-1}\phi _\lambda (v_\lambda )$.
\end{cor}

\begin{rem}
The constants $(l_1,\cdots ,l_N)$ can be computed in different ways. One can use 
\cite[\S 4.1]{caldero2} or one can compute directly from the above formula. Indeed, one knows, \cite[\S 28.1]{lusztig1}, that $(l'_1,\cdots ,l'_N):= c_\ii(v^{low}_{\lambda})$ are given by $l'_k=\langle \lambda^*, s_{i_{1}}s_{i_{2}}\cdots s_{i_{k-1}}(\alpha_{i_{k}}^{\vee})\rangle$, $1\leq k\leq N$ and resolving $(0)=b_{\ii}^{-1}\phi _\lambda c_{\ii}^{-1}(l'_1,\cdots,l'_N)$ one obtains the constants $(l_1,\cdots ,l_N)=b_{\ii}^{-1}(v^{low}_{\lambda^*} )$. These constants are given by  
$l_k=\langle \lambda, \alpha_{i_{k}}^{\vee}\rangle=\lambda_{i_{k}}$.
\end{rem}

\noindent
\textbf{Proof of the theorem:}
Fix a weight $\lambda $ in $P^+$. Let $\Phi _{\ii,\ii'}:\mathcal{C}
_\ii(\lambda )\rightarrow \Z^N$ be a family of applications indexed by two reduced words, satisfying the following properties:
\begin{enumerate}
\item[(1)] $\Phi _{\ii,\ii'}(0,\cdots,0)=b_{\ii}^{-1}\phi _\lambda
(v_\lambda )$
\item[(2)] $\Phi _{\ii,\ii'}=R_{\ii''}^\ii\circ \Phi _{\ii'',\ii'}=\Phi
_{\ii, \ii''}\circ R_{-\ii'}^{-\ii''}$
\item[(3)] If $\Phi _{\ii,\ii}(t_1,\cdots,t_N)=(t_1',\cdots ,t_N')$, then
$t'_1+t_1$ and the $t'_k$'s, $k\not=1$ only depend on $t_2, \cdots,t_N$.
\end{enumerate}

The theorem is a consequence of the following proposition:

\begin{prop}One has,
\begin{enumerate}
\item[(i)] If $(\Phi _{\ii,\ii'})$ is a family satisfying conditions (1), (2), (3), then
$$ \Phi _{\ii,\ii'}=b_{\ii}^{-1}\phi _\lambda c_{\ii'}^{-1}$$
\item [(ii)]The family $(\Phi _{\ii,\ii'})$ defined by
$$\Phi _{\ii,\ii'}(t)= [(x_\ii^{-1}\circ \zeta \circ
x_{-\ii'})^\vee]_\tr(t)+b_{\ii}^{-1}\phi _\lambda (v_\lambda )$$
satisfies conditions (1), (2), (3).
\end{enumerate}
\end{prop}

\noindent
\textbf{Proof:}
Let us first prove (ii). Given a rational subtraction-free expression $Q$, one has $[Q]_\tr(0,...,0)=0$,  so condition (1) is clear. Using \eqref{ThmTropR} and \eqref{FormZeta} the conditions (2) and (3) are also clear. 
Let us now prove (i). Let $(\Phi _{\ii,\ii'})$ be a family satisfying the conditions (1), (2), (3).
We define maps $F_\ii:\B(\lambda )\rightarrow \Z^N_{\geq 0}$ by $F_\ii(b)=\Phi _{\ii,\ii'}\circ c_{\ii'}(b)$, $b\in \B(\lambda )$. Condition (2) implies that the maps do not depend on the choice of $\ii'$.
By induction on the weight of $b$, we show that for every reduced word $\ii$ one has 
 $F_\ii(b)=b_\ii^{-1}(\phi _\lambda (b))$. Indeed, if $b=v_\lambda $, this is clear by (1).
  If $b=\tilde {f}_i(b')$, we can choose a reduced word $\ii'$ starting with $i$,
   namely $\ii'=(i,i'_2,\cdots,i_N')$. By Proposition~\ref{opeKashi}, one has 
$c_{\ii'}(\tilde {f}_i(b'))=c_{\ii'}(b')+(1,0,\cdots,0)$, so $b$ and $b'$ have the same string parameters, except the first parameters which differ by 1. Hence, by (3) we have  
$\Phi_{\ii',\ii'}\circ c_{\ii'}(\tilde {f}_i(b'))=\Phi _{\ii',\ii'}(
c_{\ii'}(b')+(1,0,\cdots,0))=\Phi _{\ii',\ii'}\circ
c_{\ii'}(b')-(1,0,\cdots,0)$ and so $F_{\ii'}(b) = F_{\ii'}(b')-(1,0,\cdots,0)$. Furthermore, by induction hypothesis one has $F_{\ii'}(b) =b_{\ii'}^{-1}(\phi
_\lambda (b'))-(1,0,\cdots,0)$. Using Proposition~\ref{opeKashi} and Theorem~\ref{baseca_et_autom}(ii), one obtains  $F_{\ii'}(b) = b_{\ii'}^{-1}(\tilde {e}_i\phi _\lambda
(b'))=b_{\ii'}^{-1}(\phi _\lambda (\tilde {f}_ib'))=b_{\ii'}^{-1}(\phi
_\lambda (b))$. Finally, by (2) one deduces $F_\ii(b)=b_\ii^{-1}(\phi _\lambda (b))$ for any reduced word $\ii$.

\subsection{Geometric lifting of $\eta_\lambda$}
Let $\ii$ be a reduced word and let $\lambda=\lambda_1\varpi_1+\cdots+\lambda_n\varpi_n$ be a dominant weight. One knows that $\ii^*$ is also a reduced word. And it is clear that the isomorphism $\delta_{\lambda^*}$ (see section \ref{definvol}) induced by $\delta$
on $V(\lambda^*)$ satisfies
$\delta_{\lambda^*} (b_\ii(t))=b_{\ii^*}(t)$.\par

Now we can give an explicit formula for the Sch\"utzenberger involution $\eta_\lambda =\delta _{\lambda ^* }\phi _\lambda $ in terms of parametrizatons of the canonical basis:

\begin{cor} \label{formuletropeta}
If $ (t_1',\cdots ,t_N')=b_{\ii^*}^{-1}\eta _\lambda
(c_\ii^{-1}(t_1,\cdots ,t_N))$, then
$$t'_k=\lambda_{i_{k}}-t_k-\sum _{j>k}a_{i_ki_j}t_j$$
\end{cor}

It is remarkable that the application 
 $b_{\ii^*}^{-1}\eta _\lambda c_\ii^{-1}$ is affine, and that its linear part does not depend on
 $\lambda $.\par 
 We define a linear map $\Omega_\ii$ as follows: for all $(\lambda,
t)=(\lambda_1,\cdots,\lambda_n,t_1,\cdots ,t_N)\in \R^{n+N}$,
$$\Omega_\ii(\lambda, t)=(\lambda_1,\cdots,\lambda_n,t_1',\cdots
,t_N'),\hspace{0.2cm} \text{where } t'_k=\lambda_{i_k}-t_k-\sum _{j>k}a_{i_ki_j}t_j$$ \label{omega}

In other words, if $b$ is an element of $\B(\lambda)$ then $\Omega_\ii(\lambda, c_{\ii}(b))=(\lambda, b_{\ii^{*}}^{-1}\eta_{\lambda}(b))$.\\

We can also give a geometric lifting of $\eta_\lambda$.\par 
Given $x \in G$, we set $\xi (x):=[\overline{w_0}(x^{-1})^{\iota}\overline{w_0^{-1}}]_{+}$.

\begin{prop}
\item[(i)] The application $\xi$ defines a bijection from $L^{w_0,e}_{>0}$ to $ L^{e,w_0}_{>0} $,
\item[(ii)] Components of $(x_{\ii}^{-1}\circ \xi \circ x_{-\ii'})^\vee$ are rational subtraction-free expressions,
\item[(iii)] One has,
$$ b_{\ii}^{-1}\eta _\lambda c_{\ii'}^{-1}(t)= [(x_\ii^{-1}\circ \xi
\circ x_{-\ii'})^\vee ]_\tr(t)+b_{\ii}^{-1}\eta _\lambda (v_\lambda ).$$
\end{prop}

\noindent
\textbf{Proof:}
It suffices to notice that $\overline{w_0}y_i(t)\overline{w_0}^{-1}=x_{i^*}(-t)$ and $\overline{w_0}t^{\alpha_{i}^{\vee}}\overline{w_0}^{-1}=t^{-\alpha_{i^*}^{\vee}}$. Then one has 
$\xi(x_{-\ii}(t_1, \cdots, t_N))=\zeta (x_{-\ii^*}(t_1, \cdots, t_N))$. The proposition then follows from Proposition \ref{FormZeta} and Theorem \ref{ThTropZeta}.

\begin{ex}
In the $\texttt{A}_2$ case, let $b$ be an element of $\B(\lambda_{1}\varpi_{1}+\lambda_{2}\varpi_{2})$. If we denote $(t_1,t_2,t_3)=c_{121}(b)$ and $\eta_{\lambda}(b)=b_{212}(t'_1,t'_2,t'_3)$ then,
\begin{equation}\left\lbrace \begin{array}{lcl}
t'_1& = & \lambda_{1}-t_1+t_2-2t_3\\
t'_2& = & \lambda_{2}-t_2+t_3\\
t'_3& = & \lambda_{1}-t_3 \end{array}
\right.
\end{equation} 
And if $(t_1,t_2,t_3)=c_{121}(b)$ and $\eta_{\lambda}(b)=b_{121}(t''_1,t''_2,t''_3)$ then,
\begin{equation}\left\lbrace \begin{array}{lcl}
t''_1& = & \lambda_{2}+t_1-t_2+2t_3-\min{(t_1+t_3,t_2)}\\
t''_2& = & \lambda_{1}-t_1-2t_3+\min{(t_1+t_3,t_2)}\\
t''_3& = & \lambda_{2}+t_3-\min{(t_1+t_3,t_2)}\end{array}
\right.
\end{equation} 
\end{ex}

\begin{ex}
In the $\texttt{B}_2$ case, let $b$ be an element of $\B(\lambda_{1}\varpi_{1}+\lambda_{2}\varpi_{2})$. If we denote $(t_1,t_2,t_3,t_4)=c_{1212}(b)$ and $\eta_{\lambda}(b)=b_{1212}(t'_1,t'_2,t'_3,t'_4)$ then,
\begin{equation}\left \lbrace \begin{array}{lcl}
t'_1& = & \lambda_{1}-t_1+t_2-2t_3+t_4\\
t'_2& = & \lambda_{2}-t_2+2t_3-2t_4\\
t'_3& = & \lambda_{1}-t_3+t_4 \\
t'_4& = & \lambda_{2}-t_4 \end{array} 
\right.
\end{equation} 

\end{ex}

\subsection{Formulas of inverse maps} 
In this section we study the inverse maps
$c_\ii\phi_\lambda b_\ii$ and $\zeta ^{-1}$. 
We have explicit formulas:

\begin{prop} \cite{CMMG} Let $\ii$ be a reduced word and let $\lambda=\lambda_1\varpi_1+\cdots+\lambda_n\varpi_n$ be a dominant weight. One has:
\begin{enumerate}
\item[(i)] If $(t_1,\cdots ,t_N)=x_{-\ii}^{-1} \zeta^{-1}x_{\ii}(t_1',\cdots
,t_N')$, then
\begin{equation} \label{inverseformula1}
t_k={t'}_k^{-1}\prod _{j>k}{t'}_j^{-a'_{i_ji_k}} ,\hspace{0.5cm} 1\leq k\leq N
\end{equation}
where $a'_{i_ji_k}=\langle \beta_{\ii,k}, \beta_{\ii,j}^\vee\rangle$.
\item[(ii)] If $(t_1,\cdots ,t_N)=c_\ii\phi_\lambda b_\ii(t_1',\cdots ,t_N')$
then
\begin{equation} \label{inverseformula2}
t_k=l'_k-t'_k-\sum _{j>k}a'_{i_ki_j}t'_j,\hspace{0.5cm} 1\leq k \leq N
\end{equation}
where $a'_{i_ji_k}=\langle \beta_{\ii,k}, \beta_{\ii,j}^\vee\rangle$ et
$l'_k=\langle\lambda ,\beta_{\ii,k}^\vee\rangle$.
\end{enumerate}
\end{prop}

Now we determine the inverse map $\zeta ^{-1}$ viewed as a map from $L^{e,w_0}_{>0}$ to $L^{w_0,e}_{>0}$
\begin{prop}
The map $\zeta ^{-1}:L^{e,w_0}_{>0}\rightarrow
L^{w_0,e}_{>0}$ is given by:
$$\zeta ^{-1}(x)=[\overline{w_0}x^T]_0x^{\iota T}$$
\end{prop}

\noindent
\textbf{Proof:}
Let us check that $\zeta ^{-1}(x)$ is well defined. We use definitions \eqref{def2cellules} and relations \eqref{relcommut}. Given $x\in L^{e,w_0}_{>0}$, set $y:=[\overline{w_0}x^T]_0x^{\iota T}$ and let us show that $y\in L^{w_0,e}_{>0}$. One can express $x$ as:
$$x=n_1h_1\overline{w_0}h_2n_2, \hspace{0.5cm} \text{with } n_1,n_2\in
N^-, h_1,h_2\in T.$$
Thus,
\begin{eqnarray*}
\overline{w_0}x^T&= \overline{w_0}n_2^Th_2^T\overline{w_0}^Th_1^Tn_1^T
&=\underbrace{\overline{w_0}n_2^T\overline{w_0}^{-1}}_{\in N^-}\underbrace{\overline{w_0}h_2^T
\overline{w_0}^{-1}h_1^T}_{\in T}\underbrace{n_1^T}_{\in N}.
\end{eqnarray*}
Hence,
\begin{eqnarray*}
y=[\overline{w_0}x^T]_0x^{\iota T}=\underbrace{\overline{w_0}h_2^T\overline{w_0}^{-1}h_1^T}_{\in T} \underbrace{n_1^{\iota T}}_{\in N}h_1^{\iota T}\overline{w_0}^{\iota T}h_2^{\iota T}n_2^{\iota T} \vspace{0.5cm}
\end{eqnarray*}

Let us commute the above two elements in $T$ and $N$:\vspace{0.5cm}
\begin{eqnarray*}
y &=&{n'}_1^{\iota T}\overline{w_0}\underbrace{h_2^T\overline{w_0}^{-1}h_1^Th_1^{\iota T}\overline{w_0}h_2^{\iota T}}_{=1}n_2^{\iota T}\\
&=&\underbrace{{n'}_1^{\iota T}}_{\in N}\overline{w_0}\underbrace{n_2^{\iota T}}_{\in N}
\end{eqnarray*}

Since $x$ is in $ N$, then $y=[\overline{w_0}x^T]_0x^{\iota T}$ is in $B^-$. Hence
$y$ belongs to $ L^{w_0,e}$ by definitions (\ref{def2cellules}). It remains to check the positivity of $y$. Fix a reduced word~$\ii$ and let us use Theorem \ref{thm.param.cells} for $x$ and $y$. One can write
 $y= x_{-\ii}(t_1,\cdots ,t_N)=y_{i_1}(t_1)t_1^{-\alpha _{i_1} ^{\vee}} \cdots
y_{i_N}(t_N)t_N^{-\alpha _{i_N} ^\vee}$ with $(t_1,\cdots ,t_N) \in \C_{\not = 0}^N$. Using relations    \eqref{relcommut} we write $y=hy_{i_1}(t'_1)\cdots y_{i_N}(t'_N)$ with a certain $h\in T$ and $t'_k=\prod_{j>k}t_j^{a_{i_ji_k}}$. But there also exists $(t''_1,\cdots ,t''_N)\in \R_{>0}^N$ such that $x=x_\ii (t''_1,\cdots,t''_N)$. Then $y=[\overline{w_0}x^T]_0x^{\iota T}= [\overline{w_0}x^T]_0y_{i_1}(t''_1)\cdots y_{i_N}(t''_N)$. Identifying both expressions one deduces that $(t_1,\cdots ,t_N)$ belongs to $\R_{>0}^N$ and thus
$y=[\overline{w_0}x^T]_0x^{\iota T}\in L^{w_0,e}_{>0}$. Hence the definition makes sense. Furthermore it is easy to see from the definitions that we have $\zeta ([\overline{w_0}x^T]_0x^{\iota T})=x$. It means that the map $x\mapsto [\overline{w_0}x^T]_0x^{\iota T}$ is precisely the inverse map of $\zeta$.

\section{Toric and semitoric degenerations of Richardson varieties}

\subsection{Problem of toric degenerations}
A complex projective variety $X$ degenerates into a toric variety, resp. semi-toric variety (i.e a variety whose irreducible components are toric varieties), if there exists a variety $\mathcal{X}$ and a regular map $\pi: \mathcal{X}\rightarrow \C$ such that  
 $\pi^{-1}(z)\cong X$ for all $z \in \C^{*}$ and $\pi^{-1}(0)= X_0$, where $X_0$ is a toric variety,
 resp. semi-toric variety. \par
 One can construct toric degenerations of varieties in the following way. 
 Denote by $R$ the algebra of regular functions on $X$. Endow $R$ with an increasing
 filtration $(R_n)_{n\geq 0}$ such that the associated graded algebra 
 $\Gr R$ is isomorphic to one algebra of semigroup. Consider $\mathcal{R}:=\oplus _{n\geq 0}R_nt^n\subset R[t]$ where $t$ is an indeterminate over $R$. One has $\mathcal{R}/{t\mathcal{R}}\simeq \oplus_{n\geq 0}
 R_{n+1}/R_n\simeq \Gr R$ and $\mathcal{R}/{(t-z)\mathcal{R}}\simeq R$, for all $z \in \C^{*}$. 
 Thus if $\mathcal{X}:=\text{Proj }\mathcal{R}$, $\pi := t$ and $X_0:= \text{Proj } \Gr R$, then we have a (flat) toric degeneration of $X$ in $X_0$.
 
\subsection{Flag varieties } \label{SectionG/B}
Consider the flag variety $G/B$ and denote by $R$ its algebra of 
homogenous coordinates. Let $\lambda _0 $ be a regular dominant weight and denote by $\f_{\lambda _{0}}$ the corresponding ample line bundle on $G/B$. Recall that:
$$R=\bigoplus_{n\in  \N} H^0(G/B, \f_{\lambda _{0}}^{\otimes n})=\bigoplus _{\lambda \in \N \cdot \lambda _0}V^* (\lambda )\otimes v_\lambda$$
where $\N . \lambda _0$ is the cone of all the multiples of $\lambda_0$.\\

From now on, we fix a regular dominant weight $\lambda_0=\lambda_1\varpi_1+\lambda_2\varpi_2+\cdots+\lambda_n \varpi_n$. One can identify each element of $\N . \lambda _0$ with the $n$-tuple of its coordinates on the fundamental weights $\varpi_1$,..., $\varpi_n$. The set $\N . \lambda _0$ is naturally identified with $\N^n$.\\

The set $\{(bv_\lambda)^* \otimes v_\lambda, b\in \B(\lambda),\lambda
\in \N . \lambda _0\}$ is the canonical basis of $R$. The product on $R$ is given 
by:
$$((bv_\lambda)^* \otimes v_\lambda)((b'v_\mu)^* \otimes
v_\mu)=(bb'v_{\lambda+\mu})^* \otimes v_{\lambda+\mu}$$

Given a reduced word $\ii$, one parametrize an element $(bv_\lambda)^* 
\otimes v_\lambda $ of the 
canonical basis of $R$ by $(\lambda,c_\ii(b))\in \Z^{n+N}_{\geq 0}$ and 
one introduce the set of all parameters:
$$\Gamma_\ii:=\{(\lambda,c_\ii(b))\in \Z^{n+N}_{\geq 0}, \lambda \in 
\N . \lambda _0,
b\in \B(\lambda)\}$$

One has

\begin{thm} \cite{littelmann} \label{ThmLittelmann1}
The set $\Gamma_\ii$ is the set of all integral points in a rational 
polyhedral convex cone of $\R^{n+N}$.
\end{thm}

We will use the notation $b_{\lambda, t}$,
$(\lambda, t)\in\Gamma_\ii$, for the element $(bv_\lambda)^* \otimes 
v_\lambda \in R$ where $ b\in
\B(\lambda)$ is such that $c_\ii(b)=t$.\par

The following multiplicative property of two elements of the canonical 
basis, is due to Caldero:

\begin{prop} \cite{Philippe}
Given $(\lambda, t), (\lambda', t')$ in $\Gamma_\ii$,
$$b_{\lambda,t}b_{\lambda',t'}=b_{\lambda+\lambda',t+t'}+\sum_{s\in
\Z^{n+N}_{\geq 
0}}d^{s}_{\lambda,t,\lambda',t'}b_{\lambda+\lambda',s}$$
with $d^{s}_{\lambda,t,\lambda',t'}\neq0\Rightarrow s\prec t+t'$, 
where
$\prec$ is the usual lexicographic order of $\Z^N_{\geq 0}$.
\end{prop}
Using this property one can construct in a first step, a $\Gamma_\ii$-filtration of $R$ such that the associated graded algebra is 
isomorphic to the algebra of the semigroup $\C[\Gamma_\ii]$. In a second step, using an adapted linear form $e:\Z^{N}_{\geq 0}\rightarrow \N$ (see \cite{Philippe}) one constructs a $\N$-filtration 
$\F_\ii:=(\F_{\ii,m})_{m\in\N}$ of $R$ such that the associated graded algebra is 
isomorphic to $\C[\Gamma_\ii]$.
Denote by $\Gr R$ the graded algebra associated to this filtration and denote by
$\bar{b}_{\lambda, t}$ the image of $b_{\lambda, t}$ in $\Gr R$. The elements $\bar{b}_{\lambda, t}$ satisfy 
$\bar{b}_{\lambda, t}\bar{b}_{\lambda', t'}=\bar{b}_{\lambda+\lambda',
t+t'}$, thus one has $\Gr R=\oplus_{(\lambda,
t)\in\Gamma_\ii}\C\bar{b}_{\lambda, t}=\C[\Gamma_\ii]$.

\begin{rem}
The toric variety Proj($\C[\Gamma_\ii]$) is the same as the toric variety constructed from the convex polytope $\cc_\ii(\lambda_0):=c_\ii(\B(\lambda_{0}))=\lambda_0\times \Z^N_{\geq 0}\cap \Gamma_\ii$. By \cite{littelmann}, one knows equations of the polytope $\cc_\ii(\lambda_0)$. One has $\cc_\ii(\lambda_0)=\{ (t_1,\cdots, t_N)\in \cc_\ii \,| \;t_k= \lambda_{i_k}-\sum _{j>k}a_{i_ki_j}t_j, \; 1\leq k\leq N\}$ where $\lambda _i$ is the coordinate of $\lambda_0$ on $\varpi_i$ (note that one can also obtain these equations with the formula \eqref{maformuletrop} using the positivity of the $t'_k$'s). In the cases where $G$ is of type $\texttt{A}_n$ or $G$ arbitrary and $\ii$ is a nice decomposition of $w_0$, one also has equations for the string cone $\cc_\ii$ (see \cite[\S 3.4]{berensteinzelevinsky2}, \cite[\S 4]{littelmann}).\\

\end{rem}

\subsection{Schubert varieties} \label{SectionXw}
Consider the Bruhat cellular decompositions:
$$G/B=\displaystyle \bigcup _{w\in W}BwB/B=\displaystyle \bigcup _{\tau\in W}B^-{\tau }B/B$$
Closures $X_w:=\overline {BwB/B}$, $w\in W$, in $G/B$, are the so-called Shubert varieties. It is well known that $\dim(X_w)=\ell(w)$. Let us denote by 
$X^{\tau }:=w_0(X_\tau)=\overline {B^-{w_0\tau }B/B}$, $\tau \in W$, the oppposite Schubert varieties in
$G/B$. Recall that $G/B=X_{w_0}=X^{w_0}$.\par
 In the sequel, let us fix two elements $w$ and $\tau$ in $W$. The algebra $R_w$ associated to $X_w$ is a quotient of $R$ by a 
 certain ideal $I_w:=\bigoplus _{\lambda \in \N . \lambda _0}V_w(\lambda )^\perp\otimes v_\lambda $, where $V_w(\lambda )^\perp $ is
 the orthogonal of $V_w(\lambda )$ in $V(\lambda)^*$.
By Theorem~\ref{ThmBaseModule} the ideal $I_w$ is compatible with the canonical basis of $R$. More precisely, 
$\{(bv_{\lambda})^*\otimes v_\lambda, \lambda \in \N . \lambda _0, b\not \in \B_w(\lambda)\}$ is a basis of $I_w$.
Denote by $\pi_w$ the canonical projection of $R$ onto $R/I_w=R_w$. 
The set $\{\pi_w(bv_\lambda^* \otimes v_\lambda), b\in \B_w(\lambda), \lambda \in \N . \lambda _0\}$ is a basis of $R_w$.
Let $\ii$ be a reduced word, define
$$\Gamma^w_\ii:=\{(\lambda,c_\ii(b))\in \Z^{n+N}_{\geq 0}, \lambda \in
\N . \lambda _0, b\in \B_w(\lambda)\}$$

We will say that a reduced word $\ii=(i_1,i_2,\cdots, i_N)$ is \textit{adapted} to $w$, if $w=s_{i_1}s_{i_2}\cdots s_{i_p}$
is a reduced decomposition.

\begin{thm} \cite{littelmann} \label{ThmLittelmann2}
If the reduced word $\ii$ is adapted to $w$ then $\Gamma^w_\ii$ is a face of the cone $\Gamma_\ii$. Moreover, $\Gamma^w_\ii =\Gamma_\ii \cap (\Z^{n+p}_{\geq 0}\times \{0\}^{N-p})$ where $p=\ell(w)$.
 In particular, $\Gamma^w_\ii$ is the set of all integer points in a rational 
polyhedral convex cone of $\R^{n+N}$.
\end{thm}

Using the above filtration $\F_\ii$, one can construct a filtration
$\F_\ii^w$ of $R_w$, as follows $\F_{\ii,m}^w:=\F_{\ii,m}+I_w$, $m\in \N$.
The associated graded algebra $\Gr R_w$ is such that $\Gr R_w=\Gr R/\Gr{I_w}$. The ideal $\Gr{I_w}$ of $\Gr R$ is generated by $\{\bar{b}_{\lambda, t}, (\lambda, t)\not \in \Gamma^w_\ii \}$. \par
In the case where $\ii$ is adapted to $w$ one has $\Gr R_w=\C[\Gamma^w_\ii]$ which is the algebra of a semigroup. Therefore one obtains toric degenerations of Shubert varieties $X_w$.\par
In the case where $\ii$ is not adapted to $w$, we will see in the next paragraph that $\Gamma^w_\ii$ is a union of faces of
$\Gamma_\ii$. Thus, we will deduce semitoric degenerations of the Schubert varieties $X_w$. \par
Regarding opposite Shubert varieties we denote by $R^\tau$ the algebra associated. This algebra is obtained from $R_\tau$ by the $w_0$-action; one has
$R^\tau=w_0(R_\tau)=R/{w_0(I_\tau)}$. The ideal $I^\tau:=w_0(I_\tau)=\bigoplus _{\lambda \in \N . \lambda _0}\eta_\lambda(V_\tau(\lambda ))^\perp \otimes v_\lambda $ of $R$ has the basis 
$\{(bv_{\lambda})^*\otimes v_\lambda, \lambda \in \N . \lambda _0, b\not \in \eta_\lambda(\B_\tau(\lambda))\}$.
Denote by $\pi^\tau$ the canonical projection of $R$ onto $R/I^\tau$. 
The elements  $\{\pi^\tau(\eta_\lambda(bv_\lambda)^*\otimes v_\lambda)$, where
$\lambda \in \N . \lambda _0$, $b\in \B_{\tau}(\lambda)\}$ constitute a basis of $R^\tau$.
We set:
$$\tilde{\Gamma}^\tau_\ii:=\{(\lambda,c_\ii(b))\in \Z_{\geq 0}^{n+N},\lambda \in \N . \lambda _0, b\in \eta_\lambda(\B_{\tau}(\lambda))\}$$
We will also show that $\tilde{\Gamma}^\tau_\ii$ is a union of faces of $\Gamma_\ii$.

\subsection{Richardson varieties} \label{SectionXwtau}
 More generally, consider Richardson varieties $X_w^\tau:=X_w\cap X^\tau$,
$w,\tau\in W$. Since $G/B=X_{w_0}=X^{w_0}$ the Schubert varieties, resp opposite Schubert variety, are the particular cases corresponding to $\tau=w_0$, resp. $w=w_0$.
Recall that if $\ell(w)+\ell(\tau)<N$ then the intersection $X_w^\tau$ is empty.\par 

For the sequel we fix $w$ and $\tau$ in $W$ such that $X_w^\tau\neq \emptyset$, and an arbitrary reduced word $\ii$. 
 Denote by 
$I^\tau_w:=I_w+I^\tau$ and $R^\tau_w:=R/{I^\tau_w}$.
The ideal $I^\tau_w$ is generated by $\{(bv_{\lambda})^*\otimes v_\lambda, \lambda \in \N . \lambda _0, b\not \in
\eta_\lambda(\B_\tau(\lambda))\cap\B_w(\lambda)\}$. We set: 
 $$\Gamma^{w,\tau}_\ii:=\Gamma_\ii^w\cap \tilde{\Gamma}^\tau_\ii=\{(\lambda,c_\ii(b))\in \Z_{\geq 0}^{n+N},\lambda \in \N . \lambda _0, b\in \eta_\lambda(\B_{\tau}(\lambda))\cap \B_w(\lambda)\}$$
  Then $I^\tau_w= \left\langle b_{\lambda,t}, (\lambda,t)\not \in \Gamma^{w,\tau}_\ii\right\rangle $.
The set $\{\pi_w^\tau((bv_\lambda)^*\otimes v_\lambda)$, $\lambda \in \N . \lambda _0$, $b\in
\eta_\lambda (\B_{\tau}(\lambda))\cap \B_w(\lambda)\}$, where $\pi_w^\tau$ is the canonical projection of
$R$ onto $R/{I_w^\tau}$, is a basis of $R_w^\tau$.\par

One constructs a filtration $\F_\ii^{w,\tau}$ of $R_w^\tau$ using the previous filtration $\F_\ii$ of $R$, as follows
$\F_{\ii,m}^{w,\tau}=\F_{\ii,m}+I_w^\tau$, $m\in \N$. The associated graded algebra $\Gr R_w^\tau$ is the quotient
$\Gr R_w^\tau=\Gr R/\Gr{I_w^\tau}$. The ideal $\Gr{I_w^\tau}$ of $\Gr R$ is generated by
 $\{\bar{b}_{\lambda, t}, (\lambda, t)\not \in \Gamma^{w,\tau}_\ii \}$. \par

\begin{prop} \label{unionface}
The set $\Gamma_\ii^{w,\tau}$ is a union of faces of $\Gamma_\ii$. In particular $\Gamma_\ii^{w,\tau}$ is the set of all integer points in a rational polyhedral convex cone of $\R^{n+N}_{\geq0}$.
\end{prop}

This proposition is a consequence of the following lemma.

\begin{lem}\label{lemmepierre}
Let $\Gamma _\R$ be a rational polyhedral convex cone of $\R^{n}_{\geq0}$ and let $\Gamma '_\R$ be an arbitrary  subset of $\Gamma _\R$. Consider $\Gamma :=\Gamma _\R \cap \Z_{\geq 0}^n $ and $\Gamma ':=\Gamma '_\R
\cap \Z_{\geq 0}^n $. If $\Gamma '$ satisfies the following conditions:
\item[(a)] For all $\gamma \in \Gamma $ and $\delta \not\in \Gamma
'$, one has $\gamma +\delta \not\in \Gamma '$ 
\item[(b)] For all $\gamma \in \Gamma '$ and $m\in \N$, one has $m\gamma 
\in \Gamma '$\\
then $\Gamma '$ is a union of faces of $\Gamma _\R$ intersected with $\Z_{\geq 0}^n$.
\end{lem}

\noindent
\textbf{Proof}:
Given $\gamma \in \Gamma '$ consider a face $\Phi _\R$ of 
$\Gamma _\R$ (eventually $\Phi _\R=\Gamma_\R$) such that $\gamma $ 
belongs to the relative interior of $\Phi_\R$. 
It suffices to prove that all elements of $\Phi:=\Phi_\R \cap \Z_{\geq 0}^n$ belongs to $\Gamma '$. Assume $\delta \in \Phi$. Since $\gamma $ lies in the interior of $\Phi_\R$, there exists $m\in \N$ such that $\gamma -\frac{1}{m}\delta $ is in $\Phi_\R$. Hence $\gamma ':=m\gamma -\delta \in \Phi_\R\cap \Z_{\geq 0}^n$, and therefore $\gamma ' \in \Gamma $. 
If $\delta \not\in \Gamma '$ then by condition $(a)$ one has $\gamma '+\delta \not\in \Gamma '$, 
but one also has $\gamma '+\delta =m\gamma \in \Gamma '$ by $(b)$, this is contradiction. 
Therefore one has $\delta \in \Gamma '$. This proves Lemma \ref{lemmepierre}.\\

\noindent
\textbf{Proof of the proposition}:
Let us apply Lemma \ref{lemmepierre} with $\Gamma^{w,\tau}_\ii\subset \Gamma_\ii$. 
Taking into account that $\Gr{I_w^\tau}=\left\langle \bar{b}_{\lambda, t}, (\lambda, t)\not \in \Gamma^{w,\tau}_\ii\right\rangle $ is an ideal of $\Gr R$, one concludes that the set $\Gamma^{w,\tau}_\ii$ satisfies condition (a) of the lemma.
Condition (b) is also satisfied. Indeed, fix $(\lambda,t)\in \Gamma ^{w,\tau}_\ii$ and $m\in \N$.\par  
First of all, choose $\ii'$ adapted to $w$ and use the reparametrization
 $R_{-{\bf i}}^{-{\bf i}'}=c_{{\bf i}'}(c_\ii)^{-1}$. One has $(\lambda,t)\in \Gamma 
^{w}_\ii$, hence $(\lambda,R_{-{\bf i}}^{-{\bf i}'}t)\in \Gamma ^{w}_{\ii'}$.
On the one hand, as a consequence of Theorem~\ref{ThmTropR} the reparametrization verifies
 $R_{-{\bf i}}^{-{\bf i}'}(mt)=mR_{-{\bf i}}^{-{\bf i}'}(t)$. On the other hand $\Gamma^w_{\ii'}$ 
 is a cone by Theorem \ref{ThmLittelmann1}. Therefore $(m\lambda,R_{-{\bf i}}^{-{\bf i}'}(mt))\in \Gamma^w_{\ii'}$.
 One easily deduces that
$m(\lambda,t)=(m\lambda,mt) \in \Gamma ^{w}_\ii$. 
Next, choose $\ii''$ adapted to $\tau$ and use the reparametrization
 $R_{-{\bf i}}^{({\bf i}'')^*}=b_{({\bf i}'')^*}^{-1}\circ c_\ii^{-1}$.
Since $(\lambda,t)\in \tilde{\Gamma} ^{\tau}_\ii$ there exists $b\in \B_\tau(\lambda)$ 
such that $(\lambda,t)=(\lambda, c_\ii \eta_\lambda(b))$. 
Thus $(\lambda, R_{-{\bf i}}^{({\bf i}'')^*}t)=(\lambda, b_{(\ii'')^*}^{-1} 
\eta_\lambda(b))= \Omega_{\ii''}(\lambda, c_\ii(b))$, where $\Omega_{\ii''}$ 
is the map defined in section \ref{omega}. One has $(\lambda, 
R_{-{\bf i}}^{({\bf i}'')^*}t)\in \Omega_{\ii''} \Gamma ^{\tau}_{\ii''}$.
Since $\Omega_{\ii''}$ is linear and $\Gamma ^{\tau}_{\ii''}$ is a cone, one has
 $(m\lambda, mR_{-{\bf i}}^{({\bf i}'')^*}t)\in \Omega_{\ii''} \Gamma ^{\tau}_{\ii''}$. 
 As a consequence of Theorem \ref{ThmTropR2} one also has 
 $mR_{-{\bf i}}^{({\bf i}'')^*}t= R_{-{\bf i}}^{({\bf i}'')^*}(mt)$.
  Hence there exists $b'\in \B_\tau(\lambda)$ such that 
  $(m\lambda ,R_{-{\bf i}}^{({\bf i}'')^*}(mt))=\Omega_{\ii''}(m\lambda,c_{\ii''}(b'))$. In other words, we have 
$(m\lambda, b_{({\bf i}'')^*}^{-1}\circ c_\ii^{-1}(mt))=(m\lambda, 
b_{(\ii'')^*}^{-1} \eta_\lambda(b'))$.
 Hence $(m\lambda,mt)=(m\lambda, c_{\ii''} \eta_\lambda(b')) \in \tilde{\Gamma} ^{\tau}_\ii$.\par 
 Finally, we obtain $(m\lambda,mt) \in \Gamma ^{w,\tau}_\ii$, that shows that the hypothesis of the
 Lemma \ref{lemmepierre} are satisfied.
We deduce that $\Gamma_\ii^{w,\tau}$ is a union of faces of $\Gamma_\ii$.\\

We are now ready to prove the main result of this section.
\begin{thm} \label{thmdegenerescence}
Richardson varieties $X_w^{\tau}$ degenerate in union of irreducible toric varieties wich are given by the faces of $\Gamma_\ii^{w,\tau}$.
\end{thm}
\noindent
\textbf{Proof:}
By Proposition \ref{unionface}, there exists faces of $\Gamma_\ii$, namely $\Phi_\ii^k$, $k=1\cdots r$, such that
 $\Gamma ^{w,\tau}_\ii=\bigcup \Phi_\ii^k$.
Thus, one has $\Gr{I_w^\tau}=\langle \bar{b}_{\lambda, t}, (\lambda, t)\not \in
\cup_{k}\Phi_\ii^k \rangle = \cap_{k} \langle \bar{b}_{\lambda, t},
(\lambda, t)\not \in \Phi_\ii^k \rangle$. Spaces
$\mathcal{I}_\ii^k:=\langle \bar{b}_{\lambda, t}, (\lambda, t)\not \in
\Phi_\ii^k \rangle$ are prime ideals of $\Gr R=\langle
\bar{b}_{\lambda, t}, (\lambda, t) \in \Gamma_\ii \rangle$ since
$\Phi_\ii^k$ are faces of $\Gamma_\ii $. The algebras $\Gr
R/{\mathcal{I}_\ii^k}=\C[\Phi_\ii^k]$ give irreducible toric varieties which are the irreducible components of the variety associated to $\Gr R_w^\tau$.

\subsection{Particular case of toric degenerations}

In the case where $G=SL_n$, one can construct toric degenerations of the Richardson varieties $X_w^\tau$ for a suitable choice of $(w,\tau)$ as follows.
\begin{prop}\label{expletorique}
Let $\ii=(1,2,1,\cdots,n,n-1, \cdots, 2,1)$ be the standard reduced word, and let $w,\tau \in W$. If $\ii$ is adapted to $w$ and
$\ii^*$ is adapted to $\tau$ then $\Gamma^{w,\tau}_\ii$ is at most one face of $\Gamma_\ii$.  
\end{prop}

\noindent
\textbf{Proof: } In the case where $\ii$ is the standard reduced word, the map $R_\ii^{-\ii}$ is linear by Corollary \ref{chgtlineaire}. The map $\Omega_{\ii^{*}}$ is always linear by definition Section \ref{omega}. By Theorem \ref{ThmLittelmann2} one knows that $\Gamma_{\ii^*}^{\tau}$ is a face of $\Gamma_{\ii^*}$, thus one deduces that $\tilde{\Gamma}_{\ii}^{\tau}=(id\times R_\ii^{-\ii})\Omega_{\ii^*}(\Gamma_{\ii^*}^{\tau})$ is a face of $\Gamma_\ii$.
Hence, $\Gamma^{w,\tau}_\ii:=\Gamma_\ii^w\cap \tilde{\Gamma}^\tau_\ii$ is only one face of $\Gamma_\ii$ when it is not empty.

%%%%%%%%%%%%%%%%%%%%%%%%%%%%%%%%%%%%%%%%%%%%%%%%%%%%%%%%%%%%%%%%%%%%%%%%%%%

\subsection{Example in the $\texttt{A}_2$ case} Let us study the case where $G=SL_3$. Let us fix $\lambda _0=\varpi_1+\varpi_2$ and $\ii=(1,2,1)$.

\begin{figure}[ht]
\centerline{\epsfbox{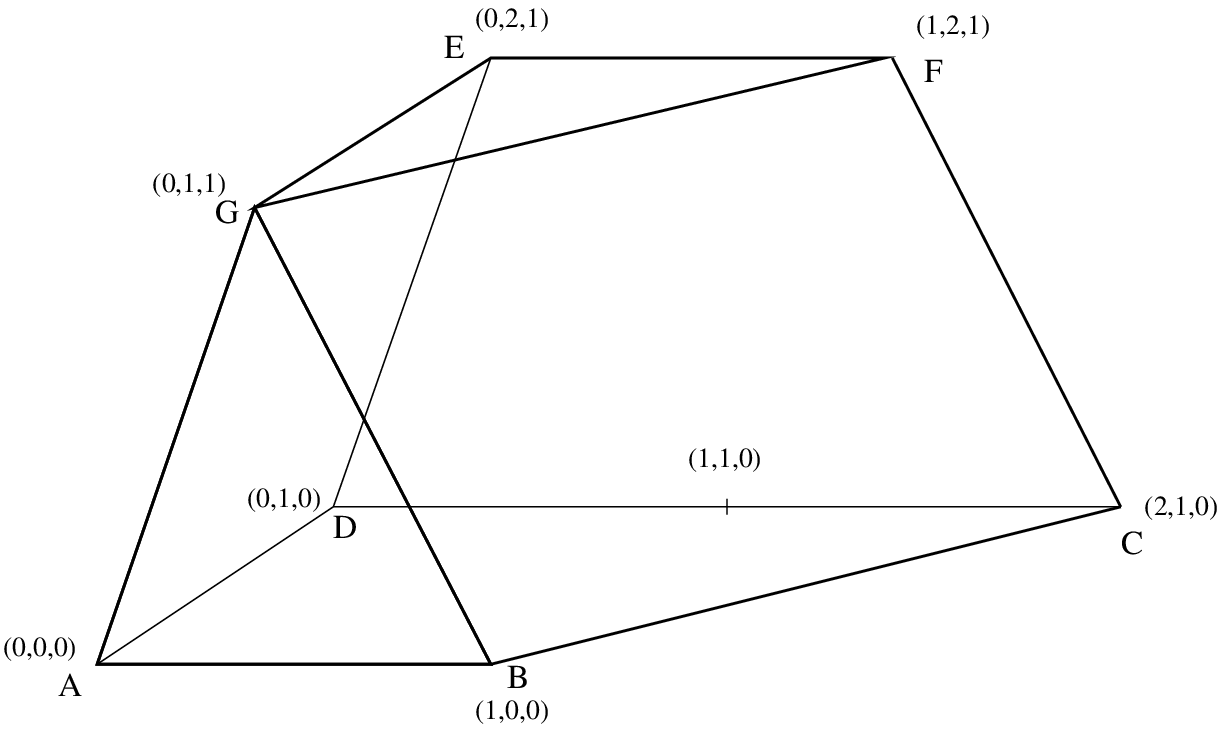}}
\caption{Polytope $c_\ii(\B(\varpi_1+\varpi_2))$}
\label{thepolytope}
\end{figure}

Figure \ref{thepolytope} represents the polytope $\cc_\ii(\lambda_0)$ of all string parameters of $\B(\varpi_1+\varpi_2)$. The cone $\Gamma_\ii$ is the cone over this polytope. The flag variety $G/B$ degenerates in the toric variety associated to this polytope.\par 
Given $w, \tau$ in $W$, denote by 
$\cc_\ii^{w,\tau}(\lambda_0):=\lambda_0\times \Z^N_{\geq 0}\cap \Gamma_\ii^{w,\tau}$ the set of all string parameters of the elements lying in $\eta_{\lambda_0}(\B_\tau({\lambda_0}))\cap\B_w({\lambda_0})$. By Theorem \ref{thmdegenerescence},  $\cc_\ii^{w,\tau}(\lambda_0)$ is a face or a union of faces (see Fig.\ref{faces}) of $\cc_\ii(\lambda_0)$ corresponding to the toric or semitoric degeneration of the subvariety $X_w^\tau$. Some of the starting varieties  $X_w^\tau$ are already toric, if they degenerate in toric varieties then the varieties are the same. \par 
Let us describe this polytope more detailly in terms of degeneration of Richardson varieties.
The vertices A, B, C, D, E, F of the polytope correspond to the $T$-fixed points $wB/B$ of $G/B$, for respectively $w=id$, $s_1$, $s_1s_2$, $s_2$, $s_2s_1$, $s_1s_2s_1$. Note that there is an extra vertex, G, resulting from the degeneration.\par
The Richardson curves are $X_{s_1}$, $X_{s_2}$, $X^{s_1}$, $X^{s_2}$, $X_{s_1s_2}^{s_1s_2}$, $X_{s_1s_2}^{s_2s_1}$, $X_{s_2s_1}^{s_1s_2}$, $X_{s_2s_1}^{s_2s_1}$, and correspond respectively to the edges [AB], [AD], [CF], [EF], [BC], [CD], [DE], [EG]~$\cup$~[BG]. The point G correspond to the intersection of the irreducible components of the degeneration of $X_{s_2s_1}^{s_2s_1}$. There are also two extra edges, namely [AG] and [FG], which will be understood as the intersection of the irreducible components of the degeneration of Richardson surfaces.\par 
The Richardson surfaces are $X_{s_1s_2}$, $X_{s_2s_1}$, $X^{s_2s_1}$, $X^{s_1s_2}$ and correspond to the faces represented Figure \ref{faces}. \par 
\medskip
\begin{figure}[!h]
%{\centerline{\input{FacesPolytope.pstex_t}}}
\begin{picture}(0,0)%
\includegraphics{PolytopeFaces.pstex}%
\end{picture}%
\setlength{\unitlength}{3108sp}%
\begingroup\makeatletter\ifx\SetFigFont\undefined%
\gdef\SetFigFont#1#2#3#4#5{%
  \reset@font\fontsize{#1}{#2pt}%
  \fontfamily{#3}\fontseries{#4}\fontshape{#5}%
  \selectfont}%
\fi\endgroup%
\begin{picture}(6678,6063)(394,-6508)
\put(631,-2671){\makebox(0,0)[lb]{\smash{{\SetFigFont{11}{13.2}{\rmdefault}{\mddefault}{\updefault}{$w=s_1s_2$}%
}}}}
\put(4006,-2626){\makebox(0,0)[lb]{\smash{{\SetFigFont{11}{13.2}{\rmdefault}{\mddefault}{\updefault}{$w=s_2s_1$}%
}}}}
\put(676,-6136){\makebox(0,0)[lb]{\smash{{\SetFigFont{11}{13.2}{\rmdefault}{\mddefault}{\updefault}{$w=w_0$}%
}}}}
\put(676,-6421){\makebox(0,0)[lb]{\smash{{\SetFigFont{11}{13.2}{\rmdefault}{\mddefault}{\updefault}{$\tau=s_2s_1$}%
}}}}
\put(3736,-6091){\makebox(0,0)[lb]{\smash{{\SetFigFont{11}{13.2}{\rmdefault}{\mddefault}{\updefault}{$w=w_0$}%
}}}}
\put(3736,-6376){\makebox(0,0)[lb]{\smash{{\SetFigFont{11}{13.2}{\rmdefault}{\mddefault}{\updefault}{$\tau=s_1s_2$}%
}}}}
\put(631,-2956){\makebox(0,0)[lb]{\smash{{\SetFigFont{11}{13.2}{\rmdefault}{\mddefault}{\updefault}{$\tau=w_0$}%
}}}}
\put(4006,-2911){\makebox(0,0)[lb]{\smash{{\SetFigFont{11}{13.2}{\rmdefault}{\mddefault}{\updefault}{$\tau=w_0$}%
}}}}
\end{picture}%
\caption{Faces $\cc_\ii^{w,\tau}(\lambda_0)$ in the polytope $\cc_\ii(\lambda_0)$}
\label{faces}
\end{figure}

The edge [AG] corresponds to the intersection of the irreducible components of $X_{s_2s_1}$, and the edge [FG] corresponds to the intersection of the irreducible components of $X^{s_1s_2}$.

\subsection{Examples in the $\texttt{B}_2$ case}

In the $\texttt{B}_2$ case there are two reduced words for $w_0$, namely $\ii=(1,2,1,2)$ and $\ii'=(2,1,2,1)$. There is no non trivial diagramm automorphism thus $\ii^*=\ii$ and ${\ii'}^*=\ii'$. We will use the word $\ii$. Using \cite[\S4]{littelmann}, one describes the string cone $\mathcal{C}_\ii$ in $Z_{\geq 0} ^4$ as follows:
\begin{equation*}
\mathcal{C}_\ii:
\begin{cases}
t_1\geq 0 & \; (\Phi_1)\\
t_2-t_3\geq 0 & \; (\Phi_2)\\
t_3-t_4\geq 0 & \; (\Phi_3)\\
t_4\geq 0 & \; (\Phi_4)
\end{cases}
\end{equation*}

Let us fix a regular dominant weight $\lambda_0=\varpi_1+\varpi_2$. By \cite[\S1]{littelmann} the polytope $\mathcal{C}_\ii(\lambda_0)$ is the intersection between $\mathcal{C}_\ii$ and the following affine cone:

\begin{equation*}
\begin{cases}
t_1-t_2+2t_3-t_4\leq 1 &\; (\tilde{\Phi}_1)\\
t_2-2t_3+2t_4\leq 1 &\; (\tilde{\Phi}_2)\\
t_3-t_4\leq 1 &\; (\tilde{\Phi}_3)\\
t_4 \leq 1 &\; (\tilde{\Phi}_4)
\end{cases}
\end{equation*}

Denote by $\Phi_i$, resp. $\tilde{\Phi}_i$, the face of the polytope $\mathcal{C}_\ii(\lambda_0)$ determined by the inequality $(\Phi_i)$, resp. $(\tilde{\Phi}_i)$.\par 
The flag variety $G/B$ degenerates in the toric variety associated to the polytope $\mathcal{C}_\ii(\lambda_0)$. The subvarieties $X^{\tau}_w$ degenerate in toric or semitoric varieties which are associated to face or union of faces of the polytope. We describe some of them.\par
The variety $X_{s_1s_2s_1}$ is associated to the face $\Phi_4$. It is a 3-dimensional polytope whose vertices are $\{(0, 0, 0, 0)$, $(0, 1, 0, 0)$, $(1, 0, 0, 0)$, $(0, 3, 1, 0)$, 
$(2, 1, 0, 0)$,  $(2, 3, 1, 0)$, $(0, 1, 1, 0)\} $. The variety $X^{s_1s_2s_1}$ is associated to the faces $\tilde{\Phi}_2\cup \tilde{\Phi}_4$. And the variety $X_{s_1s_2s_1}^{s_1s_2s_1}$  is associated to the face $\Phi_4 \cap \tilde{\Phi}_2$. It is a plane polytope whose vertices are $\{(0,1,0,0)$, $(2,1,0,0)$, $(2,3,1,0)$, $(0,3,1,0)\}$.\par
The variety $X_{s_1s_2}$ is associated to the face $\Phi_3 \cap \Phi_4$. It is a plane polytope whose vertices are $\{(0,0,0,0)$, $(0,1,0,0)$, $(2,1,0,0)$, $(1,0,0,0)\}$ (in this case the degeneration is trivial). The variety $X^{s_1s_2}$ is associated to the face $\tilde{\Phi}_2\cap \tilde{\Phi}_3$. It is a plane polytope whose vertices are $\{(1,3,2,1)$, $(2,3,1,0)$, $(0,3,1,0)$, $(0,3,2,1)\}$.\par

\end{document}